\documentclass[11pt, reqno]{amsart}

\AtBeginDocument{%
 \abovedisplayskip=11pt plus 2pt minus 8pt
 \abovedisplayshortskip=1pt plus 2pt
 \belowdisplayskip=11pt plus 2pt minus 8pt
 \belowdisplayshortskip=6pt plus 2pt minus 3pt
}

\usepackage{charter}
\usepackage{MnSymbol}
         \DeclareSymbolFont{missing}{OML}{cmr}{m}{n}
    	\DeclareMathSymbol{\ell}{\mathord}{missing}{'140}
\newcommand{\smallfrac}{\tfrac}
\newcommand{\slantfrac}{\tfrac}
\newcommand{\lrsq}{\leftrightsquigarrow}

\usepackage{thmtools}
\usepackage{pifont}
\usepackage{amsmath}
\usepackage{floatflt}
\usepackage{xcolor}
\usepackage{graphicx}
\usepackage[font = small]{caption}
\usepackage{array}
\usepackage{bbm}
\usepackage{microtype}
\usepackage{soul}
	\sodef\an{}{.06em}{.6em plus.4em}{1em plus.1em minus.1em}

\graphicspath{{./images/}}

\renewcommand{\P}{\mathbb{P}}
\newcommand{\e}{\mathrm{e}}
\newcommand{\dd}{\mathrm{d}}
\newcommand{\E}{\mathbb{E}}

\newcommand{\indi}{\mathbbm{1}}

\newcommand{\vep}{\varepsilon}
\newcommand{\Fcal}{\mathcal{F}}
\newcommand{\conn}{\longleftrightarrow}
\newcommand{\Zd}{\mathbb{Z}^d}

\newcommand{\Pp}{\P_p}

\newcommand{\Ccal}{\mathcal{C}}
\newcommand{\Tcal}{\mathcal{T}}
\newcommand{\Bcal}{\mathcal{B}}
\newcommand{\Ecal}{\mathcal{E}}
\newcommand{\Vcal}{\mathcal{V}}
\newcommand{\Xcal}{\mathcal{X}}

\newcommand{\Ppc}{\mathbb{P}_{p_c}}
\newcommand{\Epc}{\mathbb{E}_{p_c}}

\newcommand{\hits}{\longrightarrow}

\newcommand{\twa}{{(2 \wedge \alpha)}}
\newcommand{\fwa}{{(4 \wedge \alpha)}}
\newcommand{\sign}{\mathrm{sign}}

\newcommand{\lstackrel}[2]{\mathrel{\raisebox{-2.5pt}{$\stackrel{#1}{#2}$}}}
\newcommand{\ann}{\partial Q_{j, \delta r}}
\newcommand{\ck}[1]{\cite[#1]{KozNac11}}

\numberwithin{equation}{section}

\usepackage[tmargin=1.4in,bmargin=1.4in,lmargin=1.2in,rmargin=1.2in]{geometry}
\declaretheoremstyle[
    spaceabove=1em, spacebelow=1em,
    headfont=\scshape,
    notefont=\normalfont, notebraces={[}{]},
    bodyfont=\normalfont \itshape,
    postheadspace=0.5em,
    numbered=no,
                    ]{myStyle}

\theoremstyle{myStyle}

\declaretheorem[name= \underline{Theorem}]{thm}
\declaretheorem[name= \underline{Theorem}, numbered=no]{thm2}
\declaretheorem[name= Remark, style=remark]{rk}
\declaretheorem[name= \underline{Definition}]{defn}
\declaretheorem[name= \underline{Lemma}, sibling= thm]{lem}
\declaretheorem[name=\underline{Proposition}, sibling=thm]{prop}
\declaretheorem[name= \underline{Claim}, sibling= thm]{claim}

\newtheorem{innercustomthm}{\underline{Theorem}}
\newenvironment{customthm}[1]
  {\renewcommand\theinnercustomthm{#1}\innercustomthm}
  {\endinnercustomthm}

\newtheorem{innercustomlem}{\underline{Lemma}}
\newenvironment{customlem}[1]
  {\renewcommand\theinnercustomlem{#1}\innercustomlem}
  {\endinnercustomlem}

\newtheorem{innercustomdef}{\underline{Definition}}
\newenvironment{customdef}[1]
    {\renewcommand\theinnercustomdef{#1}\innercustomdef}
    {\endinnercustomdef} 

\newtheorem{innercustomclaim}{\underline{Claim}}
\newenvironment{customclaim}[1]
    {\renewcommand\theinnercustomclaim{#1}\innercustomclaim}
    {\endinnercustomclaim}

\makeatletter

\def\section{\@startsection{section}{1}%
  \z@{\linespacing\@plus\linespacing}{.5\linespacing}%
  {\normalfont\scshape}}
  
\def\subsection{\@startsection{subsection}{2}%
  \z@{.5\linespacing\@plus.7\linespacing}{-.5em}%
  {\normalfont\itshape}}

\makeatother

\usepackage[
bookmarks, colorlinks, breaklinks, pdftitle={The one-arm exponent for mean-field long-range percolation},pdfauthor={Tim Hulshof}]{hyperref}  
\definecolor{webbrown}{rgb}{.6,0,0}
\definecolor{Maroon}{cmyk}{0, 0.87, 0.68, 0.32}
\definecolor{denim}{RGB}{21, 96, 189}

\hypersetup{
    colorlinks=true, linktocpage=true, pdfstartpage=1, pdfstartview=FitV,%
    breaklinks=true, pdfpagemode=UseNone, pageanchor=true, pdfpagemode=UseOutlines,%
    plainpages=false, bookmarksnumbered, bookmarksopen=true, bookmarksopenlevel=1,%
    hypertexnames=true, pdfhighlight=/O,%
    urlcolor=webbrown, linkcolor=denim, citecolor=Maroon, %
}

\usepackage[textsize=small]{todonotes}

\begin{document}
\medskip
\begin{center}
{\sc \Large \an{The one-arm exponent for mean-field}\\
 \an{long-range percolation}}
\end{center}

\begin{abstract}
Consider a long-range percolation model on $\Zd$ where the probability that an edge $\{x,y\} \in \Zd \times \Zd$ is open is proportional to $\|x-y\|_2^{-d-\alpha}$ for some $\alpha >0$ and where $d > 3 \min\{2,\alpha\}$. We prove that in this case the one-arm exponent equals $\smallfrac12 \min\{4,\alpha\}$. We also prove that the maximal displacement for critical branching random walk scales with the same exponent. 
This establishes that both models undergo a phase transition in the parameter $\alpha$ when $\alpha =4$.
\end{abstract} 

\title[The one-arm exponent for mean-field LRP]{}
\author{Tim Hulshof}
\address{Department of Mathematics, the University of British Columbia, Vancouver, Canada.}
\address{Pacific Institute for Mathematical Sciences, the University of British Columbia, Vancouver, Canada.}
\email{thulshof@math.ubc.ca}
\date{\today}
\maketitle

{\small
\noindent
{\it MSC 2010.} 60K35, 82B27, 82B43.

\noindent
{\it Key words and phrases.}
Percolation, branching random walk, mean-field behavior, critical exponent. 
}
\vspace{1em}
\hrule
\vspace{1em}

\bigskip

\section{Introduction and main result}
In this paper we study the asymptotic behavior of the one-arm event for two closely related statistical mechanical models: branching random walk (BRW) and long-range percolation (LRP). We study LRP only in the mean-field setting. That is, we require that the dimension of the space is high enough that self-interactions become negligible. 

Roughly speaking, we say that a percolation model on $\Zd$ is long-range with parameter $\alpha$ if the edges $\{x,y\}$ of the graph with vertex set $\Zd$ and edge set $\Zd \times \Zd$ are independently removed with a probability that is asymptotically proportional to $\|x-y\|_2^{-d-\alpha}$ (where $\|\cdot\|_2$ denotes the Euclidean norm on $\Zd$). We say a retained edge is `open' and a removed edge is `closed'.
We assume throughout that the edge retention probabilities are invariant under translations and rotations by $\pi/2$. 

Our main result requires that we rigorize this notion of long-range percolation and that we state some definitions and conditions, so let us first give a simple existential version.
To state the theorem we define $Q_r$ as the intersection of the cube with sides $2r+1$ and the lattice (i.e., $Q_r := \{x \in \Zd : \| x \|_\infty \le r\}$) and its complement $Q_r^c := \Zd \setminus Q_r$ (where $\|\cdot\|_\infty$ denotes the supremum norm on $\Zd$). We write $\Ppc$ for the critical percolation measure, and $\{0 \conn Q_r^c\}$ for the event that the open cluster of the origin intersects $Q_r^c$.
\begin{thm2} Given $\alpha \in (0,2)\cup(2,\infty)$ there exist critical long-range percolation models on $\Zd$ with $d > 3\min\{2,\alpha\}$ that satisfy\emph{
\[
	c r^{-\min\{4, \alpha\}/2} \le  \Ppc(0 \conn Q_r^c) \le C r^{-\min\{4, \alpha\}/2} 
\]}
for some $C \ge c > 0$.
\end{thm2}
The theorem implies that the \emph{one-arm exponent} for high-dimensional LRP equals $\smallfrac12 \min\{4, \alpha\}$, so the model undergoes a \emph{phase transition} in $\alpha$ when $\alpha = 4$. We also show below that for critical BRW with long-range jumps the analogous exponent (for the maximal displacement) also equals $\smallfrac12 \min\{4, \alpha\}$.

It should moreover be noted that while there is another (well-known) phase transition in $\alpha$ when $\alpha=2$ that is more readily apparent (for instance in \eqref{e:twopt} and \eqref{e:volball} below), this paper presents the first rigorous proof of the existence of the `$\alpha=4$' phase transition for mean-field statistical mechanical models with long-range interactions. The phase transition at $\alpha =4$ was conjectured by Heydenreich, van der Hofstad, and the author \cite{HeyHofHul14a}, based on the observation that the maximal displacement of simple random walk on large clusters of critical LRP appears to undergo a change at $\alpha=4$. Moreover, the (easier) lower bound in the above theorem has been established by the same authors in \cite{HeyHofHul14b}. This paper, then, provides an upper bound that confirms the conjecture.
\medskip

Before we give the explicit version of this theorem, let us formally introduce the models under consideration in this paper, starting with the simpler one: critical BRW.
In this paper we view branching random walk as a random embedding of a Galton-Watson tree into $\Zd$. We will restrict ourselves to \emph{critical} Galton-Watson trees, i.e., trees with a non-trivial i.i.d.\ offspring distribution $\{p_m\}_{m=0}^\infty$ such that
\[
	\sum_{m=1}^\infty m p_m = 1\qquad \text{ and }\qquad  \sigma_p^2 := \sum_{m=2}^\infty m(m-1)p_m < \infty.
\]
Writing $\xi_v$ for the number of children of the vertex $v$ in the rooted tree $T$, the probability of $T$ under the Galton-Watson process is given by
\[
	\P(T) = \prod_{v \in T} p_{\xi_v}.
\]

We say that $\phi$ is an embedding of $T$ into $\Zd$ if $\phi$ maps the root of $T$ to the origin of $\Zd$ and $\phi$ maps any other vertex $v \in T$ to a single vertex $\phi(v) \in \Zd$. We define the branching random walk measure as the measure on configurations $(T, \phi)$ such that if $w$ is a child of $v$ in $T$, then the probability that $\phi(w)=y$ if $\phi(v)=x$ is equal to $D(x,y)$, where $D$ is the \emph{random walk one-step distribution.} That is,
\begin{equation}\label{e:BRWdef}
	\P((T,\phi)) = \P(T) \prod_{(v,w) \in T} D(\phi(v) ,\phi(w)),
\end{equation}
where the product is over pairs $(v,w)$ such that $w$ is a child of $v$ in $T$.

In this paper we will consider a the following class of one-step distributions:
\begin{defn}\label{def:D}
We always assume that $D$ is invariant under lattice symmetries and rotations by $\pi/2$ and that
\begin{equation}\label{e:Ddef}
  	D(0,x) = \frac{h(x/\Lambda)}{\sum_{x\in \Zd} h(x/\Lambda)} \qquad \forall x \in \Zd,
\end{equation}
where $\Lambda \in (0,\infty)$ and $h$ is a nonnegative bounded function on $\mathbb{R}^d$ that is piecewise continuous, has the aforementioned symmetries, satisfies
\[
	\int\limits_{\mathbb{R}^d} h(x) \dd^{d} x =1,
\]
and that is either
\begin{enumerate}
	\item supported in $[-1,1]^d$,
	\item exponentially decaying (i.e., that there exist $C, \kappa >0$ such that $h(x) \le C \e^{-\kappa \|x\|_\infty}$),
	\item or else that there exists $\alpha, c_1, c_2>0$ and $k < \infty$ such that
		\begin{equation}\label{e:asymph}
   			c_1 \|x\|_2^{-d-\alpha} \le h(x) \le c_2 \|x\|_2^{-d-\alpha} \qquad \text{ for all } \|x\|_2 \ge k.
		\end{equation}
\end{enumerate}
We consider $\alpha$ as a parameter of the model. In such cases where $h$ has bounded support or decays exponentially we set $\alpha \equiv \infty$.
\end{defn}

Our focus will be on distributions that satisfy Definition \ref{def:D}(c). The canonical example of a distribution that satisfies Definition \ref{def:D}(c) is
\begin{equation}\label{e:canon}
    D(x, y) =\mathcal{N} \max\left\{\frac{\|x-y\|_2}\Lambda, 1\right\}^{-d-\alpha}
\end{equation}
where $\mathcal{N}$ is a normalizing constant. 

Our choice to set $\alpha \equiv \infty$ for models with bounded or exponentially decaying $h$ corresponds to the fact that for these models all spatial moments are finite. Indeed, with this definition, the parameter $\alpha$ indicates which spatial moments are finite, i.e.,
\[
	\sum_{x \in \Zd} \|x\|_2^q D(0,x) < \infty \qquad \text{ iff }\qquad  q  < \alpha,
\]
for any $D$ that satisfies Definition \ref{def:D}.

We say that a branching random walk $(T,\phi)$ \emph{hits} a set $A \subseteq \Zd$ if there exists at least one vertex $v \in T$ such that $\phi(v) \in A$. For the event that a BRW started at $x$ hits the set $A$ we write $\{x \hits A\}$, and if $A = \{y\}$ for some $y \in \Zd$ we simply write $\{x \hits y\}$.
\medskip

Now we introduce the long-range percolation model. We start with the complete graph on $\Zd$, i.e., the graph $(\Zd, \Zd \times \Zd)$. We choose the \emph{percolation parameter} $p \in [0, \|D\|_{\infty}^{-1}]$ and remove each edge $\{x, y\} \in \Zd \times \Zd$ independently with probability $1 - p D(x,y)$ (where $D$ is a one-step distribution), and retain it otherwise. In this paper we only consider LRP models with one-step distributions that satisfy Definition \ref{def:D}(c).

We write $\Pp$ for the LRP measure at parameter $p$, and we write $\Ccal(x)$ for the set of all vertices that can be reached from $x$ through a path of open edges. In this paper we consider LRP at the \emph{critical point} $p_c := \inf\{p : \Pp(|\Ccal(x)| =\infty) > 0\}$. This does not depend on our choice of $x$, and it is well established that $p_c$ is a non-trivial value, i.e., LRP undergoes a phase transition in $p$ at $p_c$.
Given a vertex $x$ and a vertex set $A$, we say that $x$ is connected to $A$, and write $\{x \conn A\}$ for the event $\{\Ccal(x) \cap A \neq \emptyset\}$. If $A = \{y\}$ we simply write $\{x \conn y\}$. We write $f(n) \asymp g(n)$ if there exist $C\ge c>0$ such that $c g(n) \le f(n) \le C g(n)$ when $n$ is sufficiently large.

We make one further assumption on $D$: we assume that our choice of $D$ produces the following \emph{two-point function} asymptotics:
 \begin{equation}\label{e:twopt}
 	\Ppc(x \conn y) \asymp \|x - y\|_2^{\twa -d},
\end{equation}
where we write $(a \wedge b) := \min\{a,b\}$. Chen and Sakai \cite{CheSak12} recently proved that there exists a class of LRP models for which this bound holds when $\alpha \in (0,2) \cup (2,\infty)$ and $d > 3 \twa$ and $\Lambda$ is sufficiently large. In Remark \ref{rk:twopt} below we briefly discuss this class and the assumptions.

\begin{thm}[The one-arm exponent for BRW and LRP]\label{thm:oa} (a) Consider critical branching random walk on $\Zd$ with $\alpha \in (0,2) \cup (2,\infty)$ and $d > \twa$ and with one-step distribution $D$ that satisfies Definition~\ref{def:D}. Then,
\emph{
\[
	\P(0 \hits Q_r^c) \asymp  r^{-\fwa/2}.
\]}

(b) Consider long-range percolation on $\Zd$ with $\alpha \in (0,2) \cup (2,\infty)$ and $d > 3 \twa$ and with one-step distribution $D$ that satisfies Definition~\ref{def:D}(c) and \eqref{e:twopt}. Then,
\emph{
\[
	\Ppc(0 \conn Q_r^c) \asymp  r^{-\fwa/2}.
\]}
\end{thm}

Let us make a few remarks about the theorem:
\begin{rk}
The lower bound in Theorem \ref{thm:oa}(b) has already been proved by Heydenreich, van der Hofstad, and the author \cite[Section 5]{HeyHofHul14b}. We repeat it here for completeness, but will not prove it. For part (b) we moreover do not include the cases of Definition \ref{def:D}(a) and (b). Our proof does not work for these cases, and for models in case (a) Kozma and Nachmias \cite{KozNac11} have already proved that the one-arm exponent equals $2$. In fact, our proof is based on their work, and we will use their results throughout this paper. Proving that the one-arm exponent equals $2$ for case (b) requires several modifications to the proof of Kozma and Nachmias, but none of these are interesting enough to warrant the lengthy discussion that they require.
\end{rk}

\begin{rk} We do not treat the case `$\alpha =2$'. The reason is that our proof of Theorem~\ref{thm:oa}(a) uses asymptotics of the random walk Green's function $G(x)$ (see \eqref{e:greens} below) and our proof of Theorem \ref{thm:oa}(b) uses asymptotics of the percolation two-point function \eqref{e:twopt}, and these are not known when $\alpha=2$. When $\alpha=2$ it is conjectured (see \cite{CheSak12}) that
\[
	G(x) \asymp \frac{\|x\|_2^{2-d}}{\log \|x\|_2} \qquad \text{ and } \qquad \Ppc(0 \conn x) \asymp \frac{\|x\|_2^{2-d}}{\log \|x\|_2},
\]
when $d \ge 2$ and $d \ge 6$, respectively. Assuming these relations hold, a straightforward modification of the proof of Theorem \ref{thm:oa} shows that
\begin{equation}
	\P(0 \hits Q_r^c) \asymp \frac{\log r}{r} \qquad \text{ and } \qquad \Ppc(0 \conn Q_r^c) \asymp \frac{\log r}{r},
\end{equation}
when $d \ge 2$ and $d \ge 6$, respectively. 
\end{rk}

\begin{rk} The event $\{0 \hits Q_r^c\}$ is equivalent to the event that the \emph{maximal displacement} of the BRW exceeds $r$ (i.e., the maximal distance of any particle produced by the BRW from the starting point of the BRW). 

 There are a few related results about the maximal displacement of critical BRW. The first is due to Kesten \cite{Kest95}: Consider one-dimensional critical BRW with a translation invariant one-step distribution $D(0,x)$ that  satisfies $\sum_x x D(0,x) = 0$ and $\sum_x |x|^\alpha D(0,x) <\infty$ for some $\alpha >4$ (no symmetry assumptions are made on $D$). Let $M_n$ denote the maximal displacement of a particle at level $n$ of the tree $T$. Condition $T$ on having height at least $A n$ for some large constant $A$. Then $M_n/\sqrt{n}$ converges in distribution, and there exists a constant $C>0$ such that
\[
	\P(M_n \ge z \sqrt{n} \mid \mathrm{height}(T) > An) \le C z^{-\alpha} (1+\log(z)^2) \qquad  \forall z > 0.
\]
Kesten also proves that $M_n /\sqrt{n}$ is not tight on $\{$height$(T) > An\}$ when $\alpha \le 4$. Note that even though this result is in some sense stronger than Theorem \ref{thm:oa} when $\alpha >4$, it does not imply the upper bound when $\alpha >4$.

Kesten's result can be contrasted with an earlier result by Durrett, Kesten and Waymire \cite{DurKesWay91}, where it is proved that when $ \sum_x x D(0,x)  \in (0 , \infty)$ (i.e., the BRW has non-zero drift), then roughly speaking, the maximal displacement of a critical BRW conditioned to have $n$ particles scales to the maximum of a Brownian excursion when $\alpha>2$, while it scales to an $\alpha$-dependent function when $\alpha < 2$. In particular, there appears to be no transition at $\alpha =4$ in the presence of a drift.

In a third related result, Janson and Marckert \cite{JanMar05} determine the scaling limit of the \emph{tour of the discrete snake}. Roughly speaking, this process codifies the displacements of critical BRW conditioned on offspring size. This scaling limit is different when $\alpha < 4$ and when~$\alpha \ge 4$.
\end{rk}

\begin{rk} The restriction to the lattice $\Zd$ is easily adapted to any other transitive lattice or to $\mathbb{R}^d$. (But to achieve this one must reevaluate the bounds on moments of the expected volume of the intersection of the BRW or the percolation cluster and the cube $Q_r$, such as \eqref{e:volball} below.)
\end{rk}
\begin{rk}
It is possible to extend the result of Theorem \ref{thm:oa}(a) to dimensions $1$ and $2$ for all $\alpha \in (0,2) \cup (2,\infty)$. As it stands, the proof uses the existence of the random walk Green's function, which is not defined for recurrent random walks, whence we require that $d > \twa$. But it is not hard to see that with some minor modifications to the definition \eqref{e:zcal} below it is possible to extend these results to $d=1$ and $d=2$ by instead considering a branching random walk that is killed upon exiting a large cube $Q_R$ with $R \gg r$. 
\end{rk}
\begin{rk} The proof of the upper bound uses an induction argument that is inspired by \cite{KozNac11}, where Kozma and Nachmias prove that the one-arm exponent for high-dimensional finite-range percolation equals $2$. In \cite{KozNac11} the main difficulty of the proof comes from the complicated dependence structure of percolation. 

BRW has straightforward dependencies, and as a result parts of the proof presented here are less involved for BRW. (Other parts are more involved because we consider unbounded one-step distributions, whereas \cite{KozNac11} only considers bounded steps, which make the induction easier.)

The proof of Theorem \ref{thm:oa}(b) also uses a `regularity analysis' that is similar to the one in \cite{KozNac11}. Indeed, our Claim \ref{claim} and \cite[Theorem 2]{KozNac11} are basically the same statement. Our setting is different, but not so different that their proof fails completely. In fact, we only need to modify their proof. Considering that the proof of \cite[Theorem 2]{KozNac11} takes more than twenty pages, and that it `almost' applies to our setting, we will only give an outline of the proof of Claim~\ref{claim} (and a detailed description of the modifications to the proof of \cite[Theorem 2]{KozNac11}) in the appendix. This appendix should be read as a companion to \cite{KozNac11}. That is, we do not present a completely self-contained proof: the reader will have to refer to \cite{KozNac11} for some parts of the proof of the upper bound in Theorem \ref{thm:oa}(b).
\end{rk}

\begin{rk}\label{rk:twopt}
Chen and Sakai \cite{CheSak12} prove that the two-point function bound \eqref{e:twopt} holds for models that satisfy Definition \ref{def:D}(c) if $\alpha \in (0,2) \cup (2,\infty)$ and $d > 3 \twa$, if $\Lambda$ is sufficiently large, and if the following technical bound on the discrete second derivative of the $n$-fold discrete convolution of $D$ is satisfied:
\begin{equation}\label{e:ass1}
	|D^{*(n)} (0,x) - \smallfrac12 (D^{*(n)}(x,y) + D^{*(n)}(x,-y))| \le C n L^{\twa} \frac{(\|y\|_2 \wedge \Lambda)^2}{(\|x\|_2 \wedge \Lambda)^{ d + \twa + 2}}, 
\end{equation}
for some $C > 0$ and for all $ x \in \Zd$ and all $y \in \Zd$ such that $\|y\|_2 \le \smallfrac13 \|x\|_2$. (We write $D^{*(n)}(x,y) := \sum_{x_1, \dots, x_{n-1} \in \Zd} \prod_{i=0}^{n-1} D(x_i, x_{i+1})$ with $x_0 =x$ and $x_n =y$.)
In \cite{CheSak12} it is further shown that there exists a class of models that satisfy this bound, although it should be noted that it is not known whether \eqref{e:ass1} holds for the canonical example \eqref{e:canon}. Nevertheless, it is widely conjectured that the two-point function asymptotics \eqref{e:twopt} hold for any model that satisfies Definition \ref{def:D}(c). Supporting evidence for this conjecture is that the analogous bound holds without further assumptions for BRW when $d > \twa$, see  \eqref{e:gx} below, and for percolation models with one-step distributions as in Definition \ref{def:D}(a) and (b) when $d > 6$ and $\Lambda$ is sufficiently large.
\end{rk}

\subsection{Notation and definitions}
Parts of the proofs of Theorem \ref{thm:oa}(a) and (b) are presented separately, but they are similar and use several analogous quantities. To make these similarities visible we will denote analogous quantities for both models with the same symbols. That is, the meaning of symbols will sometimes depend on the model under consideration.

\subsubsection*{General notation}
For convenience of notation we will often write
\[
	\gamma(r) := \P(0 \hits Q_r^c) \qquad \text{ for BRW,}
\]
and similarly
\[
	 \gamma(r) := \Ppc(0 \conn Q_r^c) \qquad \text{ for LRP.}
\]
We will often state bounds in terms of the generic multiplicative constants $C$ and $c$. We always assume that $0< c,C < \infty$ but otherwise the value of these constants may change from line to line, even within equations.

When dealing with numbers, we will write $a$ rather than $\lfloor a \rfloor$ or $\lceil a \rceil$ when it is clear that the statement only holds for integers (e.g.\ in summation limits). 

\subsubsection*{Notation for BRW} We write $|T|$ for the number of vertices of the tree $T$, and in general, we write $|S|$ or $\#S$ for the cardinality of a set $S$. We say that $v \in T$ is an \emph{ancestor} of $w \in T$ if the unique path from the root of $T$ to $w$ passes through the vertex $v$. We say that $w$ is a \emph{child} of $v$ if $v$ is the direct ancestor of $w$, and write $(v,w)$. We write $|v| =n$ if $v$ is a vertex at graph distance $n$ from the root of $T$. We will abuse notation a bit and let $(T,\phi)$ also denote the \emph{random} realization of the BRW. We will also not distinguish the different probability measures typographically, and just write $\P$ for whichever measure is appropriate in the given context.

\subsubsection*{Notation for LRP} We abbreviate $\Ccal(0)$ by $\Ccal$. Sometimes we abuse notation by letting $\Ccal$ denote the induced subgraph of the vertex set $\Ccal$ in the percolation configuration. We write $\{x,y\} \in \Ccal$ for the event that $\{x,y\}$ is an open edge with both end-points in $\Ccal$. Given a configuration $\omega$ and a vertex set $A$, we define the \emph{restriction of $\omega$ to $A$} as
\[
	\omega_A(\{x,y\}) := \begin{cases} \omega(\{x,y\}) & \text{ if } x,y \in A,\\
			0 & \text{ otherwise,}
		\end{cases}
\]
for all $\{x,y\} \in \Zd \times \Zd$.
That is, $\omega_A$ is the configuration $\omega$ where all edges that do not have both end-points in $A$ have been closed. Given an event $E$ we then define the event $\{E $~on~$A\}$ as the set of all configurations $\omega$ such that $\omega_A \in E$. For connection events we write $\lstackrel{A}{\conn}$ to denote that the connection occurs on $A$. We similarly define $\{E $ off $A\} := \{E $ on $A^c\}$ and $\{E $~only~on~$A\} := E \setminus \{E $ off $A\}$.

We will frequently use the \emph{van den Berg-Kesten inequality} (BK-inequality) \cite{BerKes85}, which states that two increasing events that occur on edge-disjoint sets are negatively correlated under the percolation measure. More precisely, we write $A \circ B$ if for each $\omega \in A \cap B$ there exists a set $K(\omega) \subset \Zd \times \Zd$ such that $\omega_K \in A$ and $\omega_{K^c} \in B$. If moreover $A$ and $B$ are increasing (i.e., $\P_p(A) \le \P_q (A)$ if $p < q$), then the BK-inequality is that $\P_p (A \circ B) \le \Pp(A) \Pp(B)$.

\subsection{The structure of this paper} In Section \ref{sec:lb} we prove the lower bound in Theorem \ref{thm:oa} for BRW. Then, in Section \ref{sec:bound} we state Proposition \ref{lem:bound}, and use it to prove Theorem \ref{thm:oa} for both BRW and LRP. In Section \ref{sec:brwlem} we prove Proposition \ref{lem:bound} for BRW, and in Section \ref{sec:lrplem} we prove Proposition \ref{lem:bound} for LRP subject to Claim \ref{claim}. Finally, in Appendix \ref{app} we give an outline of the proof of Claim \ref{claim}, closely following the proof of \cite[Theorem 2]{KozNac11}.

\section{Proof of the lower bound in Theorem \ref{thm:oa}(a)}\label{sec:lb}
In this section we prove the lower bound on the one-arm probability for critical BRW. 
We prove the lower bound in two steps: first we show that the bound $\gamma(r) \ge c r^{-\alpha/2}$ using a truncation argument, and then we show that $\gamma(r) \ge c r^{-2}$ using the second-moment method.

We start with the former. If a particle of the BRW at any time takes a step of length greater than $2r+1$, then the BRW must hit $Q_r^c$. Moreover, we can restrict ourselves to BRWs that have many children, so for any $n$ we bound
\[
	\gamma(r) \ge \P(|T| \ge n+1, \exists (v,w) \in T \text{ such that } \|\phi(w) - \phi(v)\|_\infty > 2r+1).
\]
Let $Y$ be a $(0,\infty)$-valued random variable whose law is the supremum norm of a single step according to the one-step distribution $D$, i.e., $Y$ is characterized by
\begin{equation}\label{e:defx}
	\P(Y \ge t) = \sum_{x \in Q_t^c} D(0,x).
\end{equation}
It follows by independence of the steps that conditionally on the event $\{|T| \ge n+1\}$, we have the bound
\[
	\gamma(r) \ge \P(|T| \ge n+1) (1- \P(Y \le 2r+1)^n),
\]
since the event that there is a large step is bounded from below by one minus the probability that the first $n$ steps of $(T,\phi)$ are not large. Since $ 1-x \le \e^{-x}$, we can further bound
\[
	(1- (1- \P(Y > 2r+1))^n) \ge 1 - \exp\left(-n\, \P(Y > 2r+1)\right).
\]
It follows from \eqref{e:Ddef} and \eqref{e:asymph} that there exists a constant $\kappa > 0$ such that for all $r$,
\[
	\P(Y > 2r+1) = \sum_{x \in Q_{2r+1}^c} D(0,x)  \ge \kappa r^{-\alpha}.
\]
Now since $1- \e^{-x} \ge \smallfrac12 x$ when $x < 1$,
\[
	1 - \exp\left(-n \,\P(Y > 2r+1)\right) \ge \frac{\kappa n}{2 r^{\alpha}}, \qquad \text{ for all } n < \frac{2 r^{\alpha}}{\kappa}.
\]
We thus have the bound
\[
	\gamma(r) \ge \frac{\kappa n}{2 r^{\alpha}} \P(|T| \ge n+1).
\]
The classical \emph{Kolmogorov's Theorem} (see e.g.\ \cite{AthNey72}) states that
\begin{equation}\label{e:kolm}
	\P(|T| \ge s) = \frac{2}{\sigma_p \sqrt{s}}(1 \pm o(1)).
\end{equation}
Setting $n = 2 r^{\alpha}/\kappa -1$, we thus obtain the bound
\begin{equation}\label{e:smalla}
	\gamma(r) \ge \P(|T| \ge 2 r^{\alpha}/\kappa)(1-o(1)) \ge  c r^{-\alpha/2}.
\end{equation}

Now we prove the bound $\gamma(r) \ge c r^{-2}$. We use the second moment method.
Let $\Vcal(S)$ denote the set of all particles of the BRW in the set $S \subseteq \Zd$, i.e.,
\begin{equation}\label{e:zcal}
	\Vcal(S) := \{v \in T : \phi(v) \in S\}.
\end{equation}
Let $k$ be some sufficiently large number to be determined later, and define the random variable
\[
	V := |\Vcal(Q_{kr} \setminus Q_r)| = |\Vcal(Q_{kr})| - |\Vcal(Q_r)|.
\]
Clearly, $V$ is supported on $\{0,1,2,\dots\}$ and $\{V > 0\} \subseteq \{0 \hits Q_r^c\}$, so that by the Paley-Zygmund inequality,
\begin{equation}\label{e:paleyzyg}
	\gamma(r) \ge \P(V > 0) \ge \frac{\E[V]^2}{\E[V^2]}.
\end{equation}

Now define the \emph{BRW two-point function} 
\[
	\tau_n(x) := \E\Big[ \sum_{v \in T: |v|=n} \indi_{\{\phi(v) =x\}}\Big].
\]
It is an easy corollary of \cite[Proposition 2.3]{Hofs06a} that $\tau_n(x) = D^{*(n)} (0,x)$.
Therefore,
\begin{equation}\label{e:greens}
	\E[|\Vcal(Q_r)|] = \sum_{x \in Q_r} \sum_{n=0}^\infty \tau_n(x) = \sum_{x \in Q_r} \sum_{n=0}^\infty D^{*(n)}(0,x) =: \sum_{x \in Q_r} G(x),
\end{equation}
where $G(x)$ is the \emph{RW Green's function} at its critical point, with step distribution $D$. For the one-step distributions $D$ considered in this paper it is a well-established fact (see e.g.\ \cite{CheSak12, Hara08, HarHofSla03}) that when $\alpha \in (0,2) \cup (2,\infty)$ and $d > \twa$,
\begin{equation}\label{e:gx}
	G(x) \asymp  \|x\|_2^{\twa-d},
\end{equation}
so it follows that
\begin{equation}\label{e:volball}
	\E[|\Vcal(Q_r)|] \asymp \sum_{x \in Q_r}  \|x\|_2^{\twa-d} \asymp r^{\twa},
\end{equation}
and therefore, if we choose $k$ sufficiently large,
\[
	\E[V] = \E[|\Vcal(Q_{kr})|] - \E[|\Vcal(Q_r)|] \ge c (kr)^{\twa} - C r^\twa = c r^\twa.
\]

To bound $\E[V^2]$ we define the \emph{three-point function} 
\[
	\tau_{n,m}(x,y) := \E\Big[ \sum_{v \in T: |v|=n} \sum_{w \in T: |w|=m} \indi_{\{\phi(v) =x\}} \indi_{\{\phi(w) =y\}}\Big].
\] 
It similarly follows from \cite[Proposition 2.3]{Hofs06a} that
\[
	\tau_{n,m}(x,y) = \sigma_p^{2} \sum_{k=0}^{(n \wedge m)} \sum_{z \in \Zd} D^{*(k)} (z) D^{*(n-k)}(x-z) D^{*(m-k)}(y-z).
\]
Again using \eqref{e:gx} to perform this sum we get when $\alpha \in (0,2) \cup (2,\infty)$ and $d > \twa$ that
\[
	\begin{split}
	\E[|\Vcal(Q_{kr})|^2] &= \sum_{x,y \in Q_{kr}} \sum_{n,m =0}^\infty \tau_{n,m}(x,y)\\
	& = \sigma_p^2 \sum_{x,y \in Q_{kr}} \sum_{z \in \Zd} G(z) G(x-z) G(y-z)\\
	& \asymp \sum_{x,y \in Q_{kr}} \sum_{z \in \Zd} \|z\|_2^{\twa-d}  \|x-z\|_2^{\twa-d}   \|y-z\|_2^{\twa-d}.
\end{split}
\]
We evaluate the sum in two parts. If $\|z\|_2 \le 10dr$, then we first sum over $x$ and $y$ and then $z$ to obtain
\[
	 \sum_{\|z\|_2 \le 10dr} \sum_{x,y \in Q_{kr}}\|z\|_2^{\twa-d}  \|x-z\|_2^{\twa-d}   \|y-z\|_2^{\twa-d} \asymp r^{3 \twa}.
\]
If, on the other hand, $\|z\|_2 > 10dr$, then, since $\|x\|_2 < \sqrt{d} r$ for all $d \ge 1$, it follows that $\smallfrac12 \|z\|_2 < \|x-z\|_2 < 2\|z\|_2$, and similarly, $\smallfrac12 \|z\|_2 < \|y-z\|_2 < 2 \|z\|_2$. Therefore,
\[
	\sum_{x,y \in Q_{kr}} \sum_{\|z\|_2 > 10dr} \|z\|_2^{\twa-d}  \|x-z\|_2^{\twa-d}   \|y-z\|_2^{\twa-d}  \asymp r^{2d} \sum_{\|z\|_2 > 10dr} \|z\|^{3\twa -3d} \asymp r^{3 \twa}.
\]
As a result 
\begin{equation}\label{e:varball}
	\E[|\Vcal(Q_{kr})|^2]  \asymp r^{3 \twa},
\end{equation}
and it follows that $\E[V^2] \le C r^{3\twa}$ as well.

Substituting the bounds on $\E[V]$ and $\E[V^2]$ into \eqref{e:paleyzyg} we obtain
\begin{equation}\label{e:biga}
	\gamma(r) \ge \frac{c r^{2\twa}}{Cr^{3 \twa}} \ge \frac{c}{r^{2}}.
\end{equation}

Taking the maximum among the bounds \eqref{e:smalla} and \eqref{e:biga} completes the proof of the lower bound of Theorem~\ref{thm:oa}(a). \qed

\section{Proof of the upper bounds in Theorem \ref{thm:oa}}\label{sec:bound}
From here on we write $\rho$ for the long-range one-arm exponent, i.e., $\rho := \smallfrac12 (4 \wedge \alpha)$.

For both BRW and LRP we will prove the upper bound with an induction argument. The following proposition provides the crucial bound for this induction:
\begin{prop}\label{lem:bound} 
Let $\alpha \in (0,2) \cup (2,\infty)$ and choose $\beta$ such that
\emph{
\begin{equation}\label{e:betadef}
	\frac{11}{10\fwa} < \beta < \begin{cases}
	 	\frac{\alpha+1}{\alpha^2} & \text{\emph{ if }} 0 < \alpha \le 4,\\
		\frac{5}{16} & \text{\emph{ if }} \quad \alpha > 4.
		\end{cases}
\end{equation} }
For BRW with $d > \twa$ that satisfies Definition \ref{def:D}, and for LRP with $d > 3 \twa$ that satisfies Definition \ref{def:D}(c) and \eqref{e:twopt}, there exists $C_1, c_2>0$ such that for all $\lambda \in (0,1]$, and $\vep_0(\lambda)$ such that for all $\vep \in (0,\vep_0)$, 
\emph{
\begin{equation}\label{e:gammabd}
	\gamma(r(1+\lambda)) \le \frac{C_1}{\sqrt{\vep} r^\rho} + \vep^{\beta \rho} r^{\rho} \gamma(r) \gamma \left(\smallfrac12 \lambda r \right) + (1-c_2)\gamma(r).
\end{equation}
}
\end{prop}
\begin{rk} The constraints on the exponent $\beta$ are of a technical nature. They summarize a number of constraints that are needed in the proofs of Theorem \ref{thm:oa} and Proposition \ref{lem:bound}. In particular, these constraints guarantee that $\beta$ satisfies $-(\beta \wedge 1/(2\rho)) < 1-2 \beta \rho < 0$ and $2 \beta \rho  > \smallfrac{11}{10}$, which are needed in the proofs of both BRW and LRP, and that $\beta$ satisfies the constraint $-1 < \beta \left(\left( \frac{\alpha}{2\alpha + 2} \wedge \smallfrac{2}{5}\right)- \rho\right) < 0$, which is needed in the proof of LRP. 
\end{rk}

We will prove this proposition for BRW in Section \ref{sec:brwlem}, and for LRP in Section \ref{sec:lrplem}.
Before we do this, let us give the proof of the upper bounds in Theorem \ref{thm:oa}.

\proof[Proof of the upper bounds in Theorem \ref{thm:oa} subject to Proposition \ref{lem:bound}]
The proof is by induction. For ease of notation, we will prove the claim for $\gamma(r(1+\lambda))$. We assume that there exists a constant $M$ (to be determined below) such that
\begin{equation}\label{e:ass}
	\gamma(s) \le \frac{M}{s^\rho} \qquad \forall s \le r.
\end{equation}
Choose $\lambda$ sufficiently small so that 
\begin{equation}
	(1+\lambda)^\rho \le 2 \qquad \text{ and } \qquad (1-c_2) (1+\lambda)^\rho \le (1-\smallfrac12 c_2).
\end{equation}
Now choose $\vep_0(\lambda)$ as in Proposition \ref{lem:bound}. Fix $M$ so large, and $ \eta \in (0, 1)$ so that
\begin{equation}\label{e:mass1}
	M^{-2\eta} \le \vep_0 (\lambda)
\end{equation}
and
\begin{equation}\label{e:mass2}
	2C_1 M^{\eta -1} + 2\left(\frac{2}{\lambda}\right)^\rho M^{1-2\eta\beta\rho} \le \smallfrac 12 c_2 .
\end{equation}
Note that such an $M$ only exists when $\beta$ satisfies $2 \eta \beta \rho> 1$. Recall the constraints on $\beta$ in \eqref{e:betadef} and observe that $2 \beta \rho \ge \smallfrac{11}{10}$ for all $\alpha >0$, so there always exists an $\eta \in (0,1)$ such that this constraint is true.

For $s= (r(1+\lambda))^\rho \le M$ the assumption \eqref{e:ass} is true vacuously.
We advance the induction by applying Proposition \ref{lem:bound} and \eqref{e:ass}---\eqref{e:mass2} with $\vep = M^{-2 \eta}$:
\[\begin{split}
	\gamma(r(1+\lambda)) & \le \frac{C_1}{\sqrt{\vep} r^\rho} + \vep^{\beta \rho} r^\rho \gamma(r) \gamma\left(\smallfrac12 \lambda r\right) + (1-c_2)\gamma(r)\\
	& \le C_1 M^{\eta}\frac{1}{r^\rho} + \left(\frac{2}{\lambda}\right)^\rho M^{2-2 \eta \beta \rho} \frac{1}{r^{\rho}} + (1-c_2)M \frac{1}{r^\rho}\\
	& \le \frac{M}{(r(1+\lambda))^\rho} \left(2C_1 M^{\eta-1} + 2 \left(\frac{2}{\lambda}\right)^\rho M^{1-2 \eta \beta \rho} + 1-\smallfrac12 c_2\right)\\
	& \le \frac{M}{(r(1+\lambda))^\rho}.
\end{split}\]
This completes the proof of the upper bounds of Theorem \ref{thm:oa}. \qed

\section{The proof of Proposition \ref{lem:bound} for BRW}\label{sec:brwlem}
Since the tree is critical it is unlikely that $|T|$ is very large, but if $|T|$ is very large, then the BRW is likely to hit $Q_r^c$, so we fix $\vep>0$ and bound
\begin{equation}\label{e:split1}
	\gamma(r(1+\lambda)) \le \P(0 \hits Q_{r(1+\lambda)}^c, |T| \le \vep r^{2 \rho}) + \P(|T| > \vep r^{2 \rho}).
\end{equation}
For the second term in \eqref{e:split1} we use Kolmogorov's Theorem \eqref{e:kolm}, from which it follows that
\begin{equation}\label{e:kolm2}
	\P(|T| \ge \vep r^{2 \rho} ) \le \frac{C}{\sqrt{\vep}r^{\rho}},
\end{equation}
so this term contributes to the first term in \eqref{e:gammabd}.

To analyze the first term of \eqref{e:split1} we split the event according to whether $(T,\phi)$ contains a step longer than $\delta r$, with $\delta := \vep^{1/(2 \rho)}$. That is, we define the event
\[
	\Ecal_{\delta r} := \{\exists (i,j) \in T \text{ such that }\|\phi(j) -\phi(i)\|_{\infty} \ge \delta r \}.
\]

If $|T| \le \vep r^{2 \rho}$, then it is not very likely that the BRW ever makes a jump longer than $\delta r$, but if it does, then it has a good chance of hitting $Q_r^c$, so we bound
\begin{equation}\label{e:split2}
	\P(0 \hits Q_{r(1+\lambda)}^c, |T| \le \vep r^{2 \rho}) \le \P(|T| \le \vep r^{2 \rho}, \Ecal_{\delta r}) + \P(0 \hits Q_{r(1+\lambda)}^c, |T| \le \vep r^{2 \rho}, \Ecal_{\delta r}^c).
\end{equation}

We will first bound the first term.
Let $Y$ again denote a random variable whose law is characterized by \eqref{e:defx}. There exists a constant $\zeta >0$ such that for all $r$,
\begin{equation}\label{e:ylaw}
	\P(Y > r) = \sum_{x \in Q_r^c} D(x) \le \zeta r^{-\alpha}.
\end{equation}
Conditioning on $|T|=n$, the event $\Ecal_{\delta r}$ is just the event that at least one of $n$ i.i.d.\ random variables with the law of $Y$ exceeds $\delta r$, so that by the union bound and Kolmogorov's Theorem \eqref{e:kolm},
\begin{equation}\label{e:bigedges}\begin{split}
	\P(|T| \le \vep r^{2 \rho}, \Ecal_{\delta r}) &= \sum_{n=1}^{\vep r^{2 \rho}} \P(\Ecal_{\delta r} |\, |T|=n) \P(|T| = n)
	 \le \P(Y > \delta r) \sum_{n=1}^{\vep r^{2 \rho}} n \P(|T| = n)\\
	& \le \frac{\zeta}{(\delta r)^{\alpha}} \sum_{n=1}^{\vep r^{2 \rho}} \P(|T| \ge n) 
	\le \frac{C}{(\delta r)^{2 \rho}} \sum_{n=1}^{\vep r^{2 \rho}} \frac{1}{\sqrt{n}} 
	\le \frac{C\sqrt{\vep}}{\delta^{2 \rho} r^{\rho}} = \frac{C}{\sqrt{\vep} r^{\rho}}.
\end{split}\end{equation}
(Here we used for the first inequality that 
\[
	\sum_{n=1}^N \P(|T| \ge n) = \sum_{n=1}^N \sum_{i=n}^\infty \P(|T|=i) \ge \sum_{n=1}^N n \P(|T|=n),
\]
and for the last equality we used that $\delta = \vep^{1/(2\rho)}$.)
The term \eqref{e:bigedges} contributes to the first term in \eqref{e:gammabd}.

Now we bound the second term in \eqref{e:split2}. Fix $L = \vep^{\beta} r$. For all integer $j$ we define the 
shell of outer radius $j$ and thickness $\delta r$ as
\[
	\ann := \{x \in \Zd \, : \, j-\delta r < \|x\|_{\infty} \le j\} = Q_j \setminus Q_{j-\delta r},
\]
and we define the set
\begin{equation}\label{e:xjdef}
	\Xcal_j := \{v \in T: \phi(v) \in \ann \text{ and all ancestors $w$ of $v$ satisfy } \phi(w) \in Q_{j-\delta r}\}
\end{equation}
and the random variable $X_j := |\Xcal_j|$.
A moment's reflection will convince the reader that if the event $\{0 \hits Q_{r(1+\lambda)}^c\} \cap \{ |T| \le \vep r^{2 \rho}\} \cap \{ \Ecal_{\delta r}^c\}$ occurs, then either
\[
	\Bcal_1 := \{\exists j \in [r(1+\smallfrac14 \lambda), r(1+\smallfrac12 \lambda)]\text{ s.t. }0< X_j < L^{\rho}\} \cap\{ 0 \hits Q^c_{r(1+\lambda)}\} \cap \Ecal_{\delta r}^c
\]
occurs, or
\[
	\Bcal_2 := \{X_j \ge L^{\rho}\, \forall j \in [r(1+\smallfrac14 \lambda), r(1+\smallfrac12 \lambda)]\} \cap \{|T| \le \vep r^{2 \rho}\} \cap \Ecal_{\delta r}^c
\]
occurs.
Note in particular that on the event $ \{0 \hits Q_{r(1+\lambda)}^c\} \cap \Ecal_{\delta r}^c$ it cannot happen that $X_j =0$ for some (or all) $j \in [r(1+\smallfrac14 \lambda), r(1+\smallfrac12 \lambda)]$  because $\{ 0 \hits Q_{r(1+\lambda)}^c\}$ ensures that there exists a connection to $Q_{r(1+\lambda)}^c$, while $\Ecal_{\delta r}^c$ ensures that this connection cannot skip more than $\delta r$ levels at a time, and if there is a particle at level $j$, then $X_k \neq 0$ for all $k \in \{j+1,\dots, j+\delta r\}$.
\medskip

We bound $\P(\Bcal_1)$ with a regeneration argument: let $j_0$ be the first $j$ such that $0<X_j\le L^{\rho}$.  
Given $(T,\phi)$, define the pair $(\Tcal_{j_0}, \phi)$ as the BRW that is indexed by the tree $\Tcal_{j_0}$ that is induced by the set of vertices
\[
	\{v \in T : 0 \hits \phi(v) \text{ and all ancestors $w$ of $v$ satisfy } \phi(w) \in Q_{j_0 - \delta r}\}.
\]
Observe that on $\Ecal_{\delta r}^c$ this vertex set is equivalent to the set 
 \[
 	\Xcal_{j_0} \cup \{w \in T: w\text{ is an ancestor of }v \in \Xcal_{j_0}\},
\]
and that $\Xcal_{j_0}$ is the set of leaves of $\Tcal_{j_0}$ in $\partial Q_{j_0 ,\delta r}$.

We condition on the realizations of $(\Tcal_{j_0},\phi)$:
\[
	\P(\Bcal_1) = \sum_{(S,\psi) \text{ admissible}} \P \left(\{0 \hits Q_{r(1+\lambda)}^c\} \cap \Ecal_{\delta r}^c \mid (\Tcal_{j_0},\phi) = (S,\psi) \right) \P((\Tcal_{j_0},\phi) = (S,\psi)),
\]
where `$(S,\psi)$ admissible' means that $\P((\Tcal_{j_0},\phi) = (S,\psi)) > 0$. Note that `$(S, \psi)$ admissible' implies that $\{0 \hits Q_r^c\}$ occurs for $(S,\psi)$, since otherwise $j_0$ would not be well-defined, and thus $(S,\psi)$ would not be admissible.

If $\{0 \hits Q_{r(1+\lambda)}^c\} \cap \Ecal_{\delta r}^c$ occurs, then there must exist some $x \in \Xcal_{j_0}$ such that  $\{x \hits Q_{r(1+\lambda)}^c\}$ occurs. The BRW coming out of any $x \in \Xcal_{j_0}$ is an independent BRW started at $x$, since the condition $\{\Tcal_{j_0} = S\}$ tells us nothing about the BRWs of the descendants of $S$ besides their initial position. Moreover, since $(S,\psi)$ is admissible, it must be the case that $\partial S$, the set of leaves of $S$ corresponding to $\Xcal_{j_0}$, satisfies $|\partial S| \le L^\rho$. Moreover, $\|\psi(v)\|_\infty \le r(1+ \smallfrac12 \lambda)$ for all $v \in \partial S$ by our choice of $j_0$. Therefore it follows that
\[\begin{split}
	 \P\left(\{0 \hits Q_{r(1+\lambda)}^c\} \cap \Ecal_{\delta r}^c \mid (\Tcal_{j_0},\phi) = (S,\psi) \right) &\le \sum_{v \in \partial S} \P(\psi(v) \hits Q_{r(1+\lambda)}^c)\\
	 	& \le L^\rho \gamma\left(\smallfrac12 \lambda r\right),
\end{split}\]
where the second inequality also uses translation invariance and monotonicity of the one-arm event, i.e., $\P(x \hits Q_{a+b}^c) \le \P(0 \hits Q_b^c)$ for all $x \in Q_a$.
As a result,
\begin{equation}\label{e:b1bd}
	\P(\Bcal_1) \le L^{\rho} \gamma \left(\smallfrac12 \lambda r \right) \sum_{S \text{ admissible}}  \P(\Tcal_{j_0} = S)  \le L^{\rho} \gamma \left(\smallfrac12 \lambda r\right) \gamma(r),
\end{equation}
where the second inequality follows since  `$S$ admissible' implies that $\{0 \hits Q_r^c\}$. Because $L = \vep^\beta r$, this gives the second term in \eqref{e:gammabd}.
\medskip

Finally we bound $\P(\Bcal_2)$. Let 
\[
	N := \smallfrac14 \lambda (\vep^{\beta}+\delta)^{-1} = \smallfrac14 \lambda (\vep^\beta + \vep^{1/(2\rho)})^{-1}.
\]
For every integer $1 \le i \le N$ let
\[
	j_i = r + \smallfrac14 \lambda r  + i( L+ \delta r) \, \in [r(1+\smallfrac14 \lambda), r(1+\smallfrac12 \lambda)].
\]
We use these numbers below to (roughly) partition $Q_{r(1+\lambda/2)} \setminus Q_{r(1+\lambda/4)}$ into $N$ non-overlapping annuli of width $L$.

\begin{figure}
	\includegraphics[width =.7\textwidth]{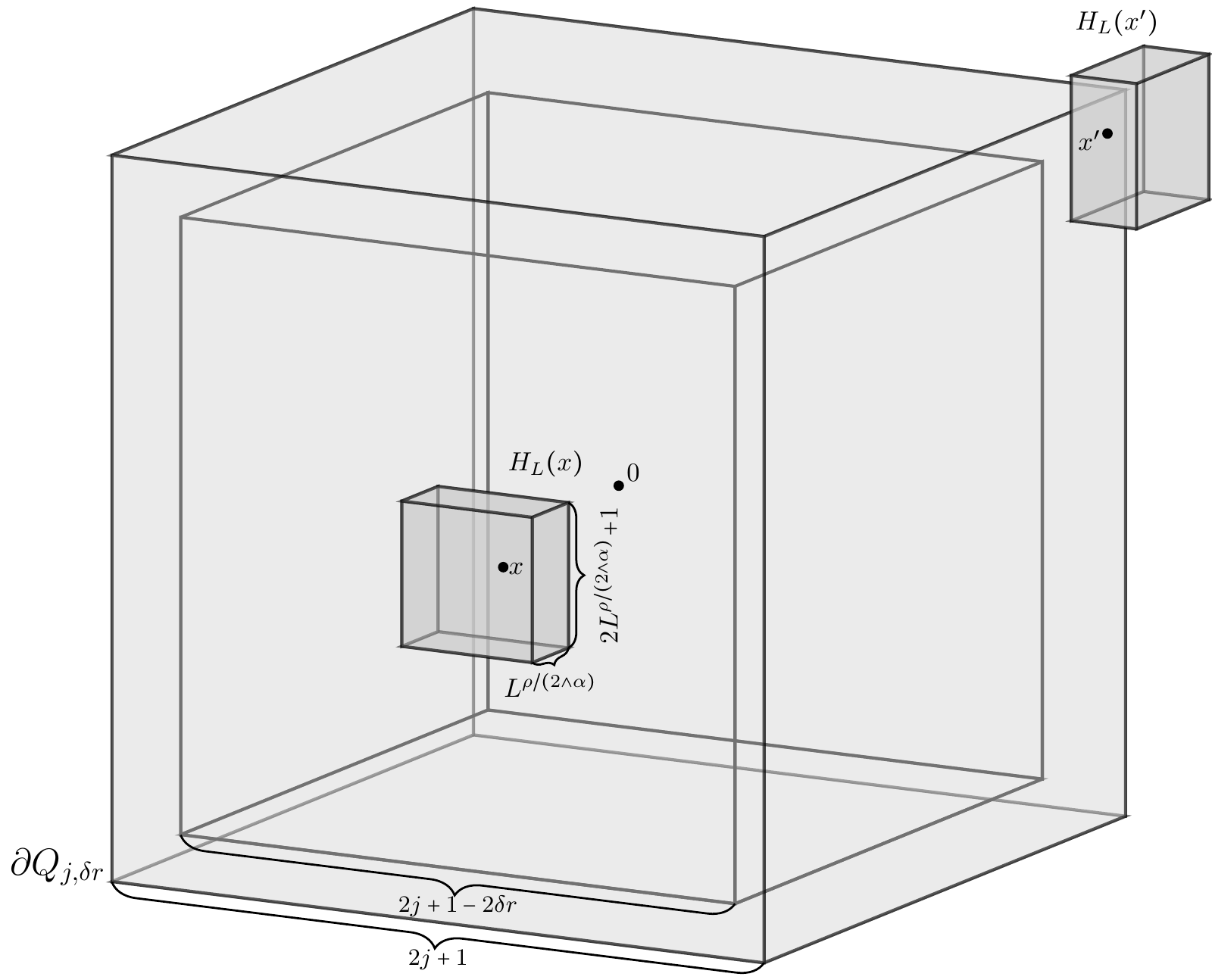}
	\caption{The annulus $\partial Q_{j,\delta r}$ with half-cubes $H_L(x)$ and $H_L(x')$ for $x, x' \in \partial Q_{j,\delta r}$.}\label{fig:hc}
\end{figure}

Let $\iota_x\in \{1, \dots d\}$ be the index of the first component of $x = (x_1, \dots, x_d) \in \Zd$ that satisfies $\|x\|_{\infty} = |x_{\iota_x}|$. Define the \emph{outward-facing half-cube}
\begin{equation}\label{e:hcube}
	H_L(x) : = \{y \in \Zd : \|x-y\|_\infty \le L^{\rho/\twa} \text{ and } |y_{\iota_x}| \ge |x_{\iota_x}| \text{ and } \sign(y_{\iota_x}) = \sign({x_{\iota_x}})\}.
\end{equation}
(See Figure \ref{fig:hc}.)

Recall the definition of $\Xcal_j$, \eqref{e:xjdef}, and define the random variables
\[
	A_j :=\# \{w \in T : \exists v \in \Xcal_j \text{ s.t.\ $v$ is an ancestor of $w$ and }\phi(w) \in H_L(\phi(v))\}
\]
and
\begin{equation}\label{e:idef}
	I := \#\{ i : X_{j_i} \ge L^{\rho} \text{ and } A_{j_i} < a L^{2 \rho}\},
\end{equation}
where $a$ is a constant that will be determined later. Observe that on $\Ecal_{\delta r}^c$, by the definition of $j_i$, $A_j$, and $H_L$, the random variables $A_{j_i}$ and $A_{j_{k}}$ count disjoint sets of particles when $i \neq k$. So if $|T| \le \vep r^{2 \rho}$, then it follows deterministically that
\begin{equation}\label{e:jdef}
	J:= \# \{i : A_{j_i} \ge a L^{2 \rho}\} < \frac{\vep r^{2 \rho}}{a L^{2 \rho}} = \frac{\vep r^{2 \rho}}{a (\vep^{\beta} r)^{2 \rho}} = a^{-1} \vep^{1-{2 \beta \rho}}.
\end{equation}
The event $\Bcal_2$ implies that all $X_{j_i} > L^{\rho}$ and that $|T| \le \vep r^{2 \rho}$. Since the total number of $j_i$'s is $N$, it follows that
\begin{equation}\label{e:nbd}
	\P(\Bcal_2) \le \P(I \ge N - J).
\end{equation}
For \eqref{e:nbd} to be non-trivial (when $\vep$ is sufficiently small) we need that $N \gg J \gg 1$, which implies that 
\begin{equation}\label{e:constraint1}
	-(\beta \wedge 1/(2\rho)) < 1-2 \beta \rho < 0.
\end{equation}
The constraints on $\beta$ in \eqref{e:betadef} ensure that this is the case.

We use Markov's inequality to bound the right-hand side of \eqref{e:nbd}. We write
\begin{equation}\label{e:Imarkov}
	\E[I] = \sum_{i=1}^{N} \P(X_{j_i} \ge L^{\rho} \text{ and } A_{j_i} < a L^{2 \rho}).
\end{equation}

We use the following claim to complete the proof:
\begin{claim}\label{claim:brw}
For BRW with $\alpha \in (0,2) \cup (2,\infty)$ and $d > \twa$ that satisfies Definition \ref{def:D}, there exists $0 < a < \slantfrac12 $ such that 
\begin{equation}\label{e:badevent}
	\P(X_{j} \ge L^{\rho} \text{\emph{ and }} A_{j} < a L^{2 \rho}) \le (1- a) \gamma(r)
\end{equation}
for all $\delta>0$ and all integer $j > r + \delta r$.
\end{claim}

Before we prove this claim, let us see how it completes the proof. By Markov's inequality, \eqref{e:Imarkov} and \eqref{e:badevent},
\begin{equation}\label{e:b2bd}
\begin{split}
	\P(\Bcal_2) &\le \P(I \ge N - J)  \le \frac{(1-a) N }{N - J } \gamma(r)\\
	& \le \frac{1-a}{1- \frac{4}{a \lambda} \vep^{1-{2 \beta \rho}}( \vep^{\beta} + \vep^{1/(2\rho)})} \gamma(r)  \le (1-c_3) \gamma(r)
\end{split}
\end{equation}
where the last inequality holds if we can choose $\vep$ small enough so that the denominator in the third inequality is arbitrarily close to one. That this is possible follows from~\eqref{e:constraint1}.

Gathering together the bounds \eqref{e:split1}, \eqref{e:kolm2}, \eqref{e:split2}, \eqref{e:bigedges}, \eqref{e:b1bd}, and \eqref{e:b2bd} gives the bound in the proposition's hypothesis. This completes the proof of Proposition \ref{lem:bound} for BRW.\qed
\medskip

\proof[Proof of Claim \ref{claim:brw}]
We write
\begin{equation}\label{e:xcond}
	\P(X_j \ge L^\rho \text{ and } A_j \le a L^{2 \rho}) = \sum_{n = L^\rho}^\infty \P(A_j \le a L^{2 \rho} \mid X_j =n) \P(X_j = n).
\end{equation}
We again use the Paley-Zygmund inequality: if $\E[A_j \mid X_j =n] > a L^{2 \rho}$, then
\begin{equation}\label{e:paleyzyg2}
	\P(A_j > a L^{2 \rho}\mid X_j =n) \ge \frac{(\E[A_j \mid X_j =n] - aL^{2\rho})^2}{\E[A_j^2 \mid X_j=n]}.
\end{equation}
Recall the definition of the half-cube $H_L$ in \eqref{e:hcube}. By independence and translation invariance of the BRW, the symmetries of $D$, and \eqref{e:volball}, there exists a constant $0 < b \le 1$ such that
\begin{equation}\label{e:zjbd}
	\E[A_j \mid X_j =n] = n\, \E[|\Vcal(H_L(0))|] \ge \smallfrac12 n  \E[|\Vcal(Q_{L^{\rho/\twa}})|] \ge bn L^\rho.
\end{equation}
Again using translation invariance and independence, the bounds \eqref{e:volball} and \eqref{e:varball}, and using in addition the fact that two vertices counted in $A_j$ can either have the same parent or different parents in $\Xcal_j$, we bound 
\begin{equation}\label{e:zj2bd}
\begin{split}
	\E[A_j^2 \mid X_j =n] &= n(n-1) \E[|\Vcal(H_L(0))|]^2 + n \E[|\Vcal(H_L(0))|^2] \\
	& \le n^2 \E[|\Vcal(Q_{L^{\rho/\twa}})|]^2 + n \E[|\Vcal(Q_{L^{\rho /\twa}})|^2]\\
	& \le C n^2 L^{2 \rho} + C' n L^{3 \rho} \le K n^2 L^{2 \rho},
\end{split}
\end{equation}
where the last bound only holds when $n \ge L^\rho$. We can assume that $K \ge 1$ and choose $a =b/(K+1)$ so that the condition $\E[A_j \mid X_j =n] > a L^{2\rho}$ is satisfied. Substituting \eqref{e:zjbd} and \eqref{e:zj2bd} into \eqref{e:paleyzyg2} we obtain that for all $n \ge L^\rho$,
\[
	\P(A_j \le a L^{2 \rho} \mid X_j =n) \le 1- \frac{(b-a) n^2 L^{2\rho}}{K n^2 L^{2 \rho}} = 1- a.
\]
Substituting this bound into \eqref{e:xcond} we get
\[
	\P(X_j \ge L^\rho \text{ and } A_j \le a L^{2 \rho}) \le (1-a)  \P(X_j \ge L^\rho) \le (1-a) \gamma(r),
\]
since $\{X_j \ge L^\rho\} \subseteq \{0 \hits Q_r^c\}$ whenever $j > r +\delta r$. This completes the proof of Claim~\ref{claim:brw}. \qed

\section{The proof of Proposition \ref{lem:bound} for LRP}\label{sec:lrplem}
The proof of Proposition \ref{lem:bound} for LRP that we give in this section is similar to the proof for BRW in the previous section, so we will supply fewer details in some of the calculations. It should however be noted that LRP is a more complicated model than BRW, because it has complicated dependencies. These dependencies will mostly manifest themselves in Claim \ref{claim} below.

To prove the one-arm exponent for high-dimensional finite-range percolation, \cite{KozNac11} uses the Barsky-Aizenman bound \cite{BarAiz91, HarSla90a}
\begin{equation}\label{p:hhs}
	\Ppc(|\Ccal| \ge n ) \asymp \frac{1}{\sqrt{n}}
\end{equation}
to establish an ersatz upper bound $\gamma(r) \le \Ppc(|\Ccal| \ge r/\Lambda) \le C/\sqrt{r}$ (where $\Lambda$ is the maximal edge length). Although \eqref{p:hhs} also holds for LRP \cite[Theorem 1.3]{HeyHofSak08}, the bound $\gamma(r) \le C/\sqrt{r}$ does not, because there is no maximal edge length (and because it contradicts the lower bound in Theorem \ref{thm:oa}(b) when $\alpha < 1$). 

In the proof of Proposition \ref{lem:bound} we will need an ersatz upper bound on $\gamma(r)$ in several places.
The following lemma establishes such a bound:
\begin{lem}\label{lem:subopt}
For LRP with $d > 3 \twa$ that satisfies Definition \ref{def:D}(c), let $\xi = \left( \frac{\alpha}{2\alpha + 2} \wedge \smallfrac{2}{5}\right)$. There exists $C>0$ such that
\begin{equation}\label{e:subopt}
	\gamma(r) \le \frac{C}{r^\xi}.
\end{equation}
\end{lem}
\proof
We use \eqref{p:hhs} to bound
\begin{equation}\label{e:subopt2}\begin{split}
	\gamma(r) &\le \Ppc(|\Ccal| > r^{2\xi}) + \Ppc(|\Ccal| \le r^{2\xi}, 0 \conn Q_r^c)\\
	& \le \frac{C}{r^{\xi}} + \Ppc(|\Ccal| \le r^{2\xi}, \exists \{v,w\} \in \Ccal \text{ such that } \|v-w\|_\infty > r^{1-2\xi}).
\end{split}
\end{equation}

To bound the second term we define the $\ell$\emph{-truncated cluster} $\widetilde\Ccal_\ell$ as the subgraph of $\Ccal$ with vertex set
\[
     \{x\in \Ccal \,:\, 0 \conn x \text{ with edges of $\| \, \cdot \,\|_\infty$-length }< \ell\},
\]
and with the set of edges that this vertex set induces in $\Ccal$ and that are shorter than $\ell$ as its edge set.
That is, $\widetilde\Ccal_{\ell}$ is the part of the cluster that can be reached from $0$ by only using open edges of length less than $\ell$. 

If $\Ccal$ contains an edge $\{v,w\}$ longer than $\ell$, then there either exist two edge-disjoint paths from $v$ to $w$, or else $|\widetilde\Ccal_\ell| < |\Ccal|$. In case of the former, there must then exist a vertex $x \in Q_{r^{\ell}/2}^c(v)$ such that $\{v \conn x\}\circ\{v \conn x\}$ occurs. So by the union bound, translation invariance, and the BK-inequality,
\begin{multline}\label{e:singledouble}
	\Ppc(|\Ccal| \le r^{2\xi}, \exists \{v,w\} \in \Ccal \text{ such that } \|v-w\|_\infty > r^{1-2\xi}) \\
	\le \Ppc(|\widetilde\Ccal_{r^{1-2\xi}}| < |\Ccal| \le r^{2\xi}) + \sum_{x \in Q_{r^{1-2\xi}/2}} \Ppc(0 \conn x)^2.
\end{multline}

To bound the second term we use that $\sum_{x} \|x\|_2^{\twa} \Ppc(0 \conn x)^2 < C$ \cite[Proposition~2.5]{HeyHofHul14b},
\begin{equation}\label{e:prop25bd}
	 \sum_{x \in Q_{r^{1-2\xi}/2}} \Ppc(0 \conn x)^2 \le C \sum_{x \in Q_{r^{1-2\xi}/2}} \frac{\|x\|_2^{\twa}}{r^{(1-2\xi)\twa}}\Ppc(0 \conn x)^2 \le \frac{C}{r^{(1-2\xi)\twa}}.
\end{equation}
(The constant in the first inequality may depend on $d$ and takes care of the discrepancy between $\| \,\cdot\,\|_2$ and $\|\,\cdot\,\|_\infty$.)
	
To bound the first term of \eqref{e:singledouble} we condition on the size of $|\widetilde\Ccal_{r^{1-2\xi}}|$. Under this conditioning, independent percolation models have a certain independence for short and long edges: Knowing $\widetilde\Ccal_{r^{1-2\xi}}$ tells us which sites can be reached from the origin via `short' edges (i.e., edges of length less than $r^{1-2\xi}$) but does not provide us with any information about the `long' edges attached to these vertices. Since moreover the status of edges is independent under $\Ppc$, and $\Ppc$ is translation invariant, we can bound
\begin{equation}\label{e:condtrunc}
\begin{split}
	 \Ppc(|\widetilde\Ccal_{r^{1-2\xi}}| < |\Ccal| \le r^{2\xi})
	 &= \sum_{n=2}^{r^{2 \xi}} \sum_{s =1}^{n-1} \Ppc(|\Ccal| =n \mid |\widetilde\Ccal_{r^{1-2\xi}}| =s ) \Ppc(|\widetilde\Ccal_{r^{1-2\xi}}| =s )\\
	 &\le \sum_{n=2}^{r^{2 \xi}} \sum_{s =1}^{n-1} s \, \Ppc( \exists v \in Q_{r^{1-2\xi}}^c \text{ s.t. } \{0,v\} \text{ is open}) \Ppc(|\widetilde\Ccal_{r^{1-2\xi}}| =s ),
\end{split}
\end{equation}
where we also used the union bound.

Using the union bound and and \eqref{e:ylaw} we get
\begin{equation}\label{e:longedgeex}
	\Ppc( \exists v \in Q_{r^{1-2\xi}}^c \text{ s.t. } \{0,v\} \text{ is open}) 
		 \le \sum_{x \in Q_{r^{1-2\xi}}^c} p_c D(0,x) \le C r^{-(1-2\xi)\alpha}.
\end{equation}

Substituting this bound into \eqref{e:condtrunc} and using \eqref{p:hhs} we get the upper bound
\begin{equation}\label{e:2ndbd}\begin{split}
	Cr^{-(1-2\xi)\alpha} \sum_{n=2}^{r^{2 \xi}} \sum_{s =1}^{n-1} s\, \Ppc(|\widetilde\Ccal_{r^{1-2\xi}}| =s )  &\le Cr^{-(1-2\xi)\alpha} \sum_{n=2}^{r^{2 \xi}} \Ppc(|\widetilde\Ccal_{r^{1-2\xi}}| \ge n )\\
		& \le Cr^{-(1-2\xi)\alpha} \sum_{n=2}^{r^{2 \xi}} \Ppc(|\Ccal | \ge n )\\
		& \le  Cr^{-(1-2\xi)\alpha} r^{ \xi} = \frac{C}{r^{\alpha - 2\alpha \xi - \xi}}.
\end{split}
\end{equation}

Finally, combining the bounds \eqref{e:subopt2}, \eqref{e:prop25bd}, and \eqref{e:2ndbd}, our choice of $\xi$ gives
\[
	\gamma(r) \le C \max \left\{r^{-\xi}, r^{-(1-2\xi)\twa}, r^{-\alpha + 2 \alpha \xi + \xi}\right\} = \frac{C}{r^\xi}.\qed
\]
\medskip

\proof[Proof of Proposition \ref{lem:bound} for LRP subject to Claim \ref{claim} below]
Recall that $\rho = \smallfrac12 \fwa$. It follows from Lemma~\ref{lem:subopt} that it suffices to show that Proposition~\ref{lem:bound} holds for $\vep \le r^{\xi - \rho} \ll 1$, since
\[
	\gamma(r(1+\lambda)) \le \frac{C}{r^{\xi}} \le \frac{C}{\sqrt{\vep}r^\rho}.
\]

From here on we will thus assume that $\vep> r^{\xi-\rho}$ (this will be important for the validity of Claim \ref{claim} below). Using \eqref{p:hhs} again we bound
\begin{equation}\label{p:split1}
	\gamma(r(1+\lambda)) \le  \frac{C}{\sqrt{\vep}r^\rho} + \Ppc(0 \conn Q_{r(1+\lambda)}^c, |\Ccal| \le \vep r^{2 \rho}),
\end{equation}
so the first term here contributes to the first term in \eqref{e:gammabd}.

To analyze the second term of \eqref{p:split1} we split the event according to whether $\Ccal$ contains an open edge longer than $\delta r$, with $\delta := \vep^{1/(2 \rho)}$. That is, we define the event
\[
	\Ecal_{\delta r} := \{\exists \{x,y\} \in \Ccal \text{ such that } \|x -y\|_{\infty} \ge \delta r \}.
\]
We bound
\begin{equation}\label{p:split2}
	\Ppc(0 \conn Q_{r(1+\lambda)}^c, |\Ccal| \le \vep r^{2 \rho})
	 \le \Ppc(|\Ccal| \le \vep r^{2 \rho}, \Ecal_{\delta r}) + \Ppc(0 \conn Q_{r(1+\lambda)}^c, |\Ccal| \le \vep r^{2 \rho}, \Ecal_{\delta r}^c).
\end{equation}

We start with the first term.
Consider the $\delta r$-truncated cluster $\widetilde\Ccal_{\delta r}$.
As explained in the proof of Lemma \ref{lem:subopt} above, the long and short edges are conditionally independent on $\widetilde\Ccal_{\delta r}$, i.e.,
\begin{equation}\label{e:condinv}
\begin{split}
        \Ppc(| \Ccal | \le \vep r^{2 \rho}, \Ecal_{\delta r}) &\le \Ppc(| \widetilde\Ccal_{\delta r} | \le \vep r^{2 \rho}, \Ecal_{\delta r})  = \sum_{n=1}^{\vep r^{2 \rho}} \Ppc(\Ecal_{\delta r} \mid |\Ccal_{\delta r}|=n) \Ppc(|\widetilde\Ccal_{\delta r}|=n)\\
       &\le \Ppc(\exists x \in Q_{\delta r}^c \text{ s.t. } \{0, x\} \text{ is open} ) \sum_{n=1}^{\vep r^{2 \rho}} n \Ppc(|\widetilde\Ccal_{\delta r}|=n)\\
       & \le \frac{C}{(\delta r)^\alpha}\, \sum_{n=1}^{\vep r^{2 \rho}}n \Ppc(|\widetilde\Ccal_{\delta r}| = n),
\end{split}
\end{equation}
where the second inequality follows from translation invariance of the model and the union bound, and the third from \eqref{e:longedgeex}. 
Applying \eqref{p:hhs} to the right-hand side of \eqref{e:condinv} again yields
\begin{equation}\label{p:bigedges}
\begin{split}
        \Ppc(| \Ccal | \le \vep r^{2 \rho}, \Ecal_{\delta r}) &\le \frac{C}{(\delta r)^{2 \rho}}\, \sum_{n=1}^{\vep r^{2 \rho}} \Ppc(|\widetilde\Ccal_{\delta r}| \ge n) \le \frac{C}{(\delta r)^{2 \rho}}\, \sum_{n=1}^{\vep r^{2 \rho}} \Ppc(|\Ccal| \ge n) \\
        & \le \frac{C}{(\delta r)^{2 \rho}}\, \sum_{n=1}^{\vep r^{2 \rho}} \frac{1}{\sqrt{n}}\le \frac{C \sqrt{\vep}}{\delta^{2 \rho} r^\rho} \le \frac{C}{\sqrt{\vep} r^{\rho}},
\end{split}\end{equation}
where in the last equality we used that $\delta = \vep^{1/(2\rho)}$.
The term \eqref{p:bigedges} contributes to the first term in \eqref{e:gammabd}.

Now we bound the second term in \eqref{p:split2}. Fix $L = \vep^{\beta} r$. Observe that the constraints \eqref{e:betadef} on $\beta$ ensure that $1 \ll L \ll r$ when $\vep \ge r^{\xi-\rho}$ (that $L \ll r$ simply follows from $\beta>0$). We define the set
\begin{equation}\label{p:xjdef}
	\Xcal_j := \left\{x \in \ann\,:\, 0 \lstackrel{Q_{j}}{\conn} x \right\}
\end{equation}
and the random variable $X_j := |\Xcal_j|$.
If the event $\{0 \conn Q_{r(1+\lambda)}^c\} \cap \{ |\Ccal| \le \vep r^{2 \rho}\} \cap \Ecal_{\delta r}^c$ occurs, then either
\[
	\Bcal_1 := \{\exists j \in [r(1+\smallfrac14 \lambda), r(1+\smallfrac12 \lambda)]\text{ s.t. }0< X_j < L^{\rho}\} \cap\{ 0 \conn Q_{r(1+\lambda)}^c\} \cap \Ecal_{\delta r}^c
\]
occurs, or
\[
	\Bcal_2 := \{X_j \ge L^{\rho}\, \forall j \in [r(1+\smallfrac14 \lambda), r(1+\smallfrac12 \lambda)]\} \cap \{|\Ccal| \le \vep r^{2 \rho}\} \cap \Ecal_{\delta r}^c
\]
occurs.

We start with the bound on $\Ppc(\Bcal_1)$: let $j_0$ be the first $j$ such that $0<X_j\le L^{\rho}$. 
As in \cite{KozNac11}, we define the \emph{cluster until level $j$} as $\Ccal_j : = \Ccal(0; Q_{j})$ as the connected cluster of the origin in the modified configuration where all edges that do not have both vertices in $Q_j$ are closed. That is, an edge $\{x,y\}$ is in $\Ccal_j$ if $\{x,y\}$ is in the cluster of the origin in $\omega_{Q_j}$. 
Observe that $\Xcal_{j_0}$ is the set of vertices of $\Ccal_{j_0}$ in $\partial Q_{j_0, \delta r}$.

We condition on the realizations of $\Ccal_{j_0}$ (i.e., we condition on all the information needed to determine $\Ccal_{j_0}$: the cluster until level $j_0$, and all the closed edges touching this cluster), to obtain
\[
	\Ppc(\Bcal_1) = \sum_{S \text{ admissible}} \Ppc \left(\{0 \conn Q_{r(1+\lambda)}^c \} \cap \Ecal_{\delta r}^c \mid \Ccal_{j_0} = S \right) \Ppc(\Ccal_{j_0} = S),
\]
where `$S$ admissible' means that $\Ppc(\Ccal_{j_0}= S) > 0$. Observe that by the definition of $\Ccal_{j_0}$, any admissible $S$ must satisfy $S \cap Q_r^c \neq \emptyset$. If $\{0 \conn Q_{r(1+\lambda)}^c\} \cap \Ecal_{\delta r}^c$ occurs, then there must exist a connection from $x$  to $Q_{r(1+\lambda)}^c$ `off $\Ccal_{j_0}$'  for some $x \in \Xcal_{j_0}$, so if we write $\partial S = S \cap \partial Q_{j_0, \delta r}$, then
\[
	\Ppc \left(\{0 \conn Q_{r(1+\lambda)}^c\} \cap \Ecal_{\delta r}^c \mid \Ccal_{j_0} = S \right) \le \sum_{x \in \partial S} \Ppc \left(x \conn Q_{r(1+\lambda)}^c \text{ off }S \mid  \Ccal_{j_0} = S \right).
\]
We can tell whether $j = j_0$ by inspecting $\Ccal_j$ (since $j_0$ is the first $j$ such that $0 < X_j < L^\rho$), so conditioning on $\{\Ccal_{j_0} = S\}$ reveals no information about the configuration `off $S$'. Therefore,
\[
	\Ppc \left(x \conn Q_{r(1+\lambda)}^c \text{ off }S \mid  \Ccal_{j_0} = S \right) = \Ppc \left(x \conn Q_{r(1+\lambda)}^c \text{ off }S \right)
	\le \Ppc \left(x \conn Q_{r(1+\lambda)}^c  \right),
\]
where the inequality follows by monotonicity of the one-arm event.
It follows that
\[
	 \Ppc\left(0 \conn Q_{r(1+\lambda)}^c \mid \Ccal_{j_0} = S \right) \le \sum_{x \in \partial S} \Ppc \left(x \conn Q_{r(1+\lambda)}^c \right) \le L^\rho \gamma\left(\smallfrac12 \lambda r\right),
\]
where the second inequality also uses translation invariance and monotonicity of $\gamma(r)$. As a result,
\begin{equation}\label{p:b1bd}
	\Ppc(\Bcal_1) \le L^{\rho} \gamma \left(\smallfrac12 \lambda r \right) \sum_{S \text{ admissible}}  \Ppc(\Ccal_{j_0} = S)  \le L^{\rho} \gamma \left(\smallfrac12 \lambda r\right) \gamma(r),
\end{equation}
since `$S$ admissible' implies that $\{0 \conn Q_r^c\}$. This gives the second term in \eqref{e:gammabd}.
\medskip

Finally we bound $\Ppc(\Bcal_2)$. The same computations as in Section \ref{sec:brwlem} will complete the proof, but we do need a percolation variant of Claim \ref{claim:brw}.
Let 
\[
	N := \smallfrac14 \lambda (\vep^{\beta}+\delta)^{-1} = \smallfrac14 \lambda (\vep^{\beta}+\vep^{1/(2\rho)})^{-1}.
\]
Recall that $L = \vep^\beta r$ and for every integer $1 \le i \le N$ let
\[
	j_i = r + \smallfrac14 \lambda r  + i( L+ \delta r) \, \in [r(1+\smallfrac14 \lambda), r(1+\smallfrac12 \lambda)].
\]
Recall the definition of $\Xcal_j$, \eqref{p:xjdef}, and define the random variable
\[
	A_j :=\# \{y \in Q_{j+L} \setminus Q_{j} : 0 \conn y\}.
\] 

\begin{claim}\label{claim}
For LRP with $\alpha \in (0,2) \cup (2,\infty)$ and $d > 3 \twa$ that satisfies Definition \ref{def:D}(c) and \eqref{e:twopt} there exist constants $0 < a < \slantfrac12$  such that for any $j$ sufficiently large, any $\theta \in (0,1)$, and any $L \ge j^{\theta}$,
\begin{equation}\label{p:badevent}
	\Ppc(X_j \ge L^{\rho} \text{\emph{ and }} A_j < a L^{2 \rho}) \le (1- a) \gamma(r).
\end{equation}
\end{claim}

With this claim we can complete the proof of Proposition \ref{lem:bound} for LRP. Earlier, we chose $\vep > r^{\xi -\rho}$, $L= \vep^\beta r$, and $j \in [r(1+\smallfrac14 \lambda), r(1+\smallfrac12 \lambda)]$, so by our constraints on $\beta$ in \eqref{e:betadef} we can find a $\theta>0$ such that $L \ge j^\theta$. Define $I$ and $J$ as in the BRW case (see \eqref{e:idef} and \eqref{e:jdef}). Then, by Markov's inequality and Claim \ref{claim},
\begin{equation}\label{p:b2bd}
\begin{split}
	\Ppc(\Bcal_2) &\le \Ppc(I \ge N - J)  \le \frac{(1-a) N }{N - J } \gamma(r)\\
	& \le \frac{1-a}{1- \frac{4}{a \lambda} \vep^{1-{2 \beta \rho}}( \vep^{\beta} + \vep^{1/(2\rho)})} \gamma(r)  \le (1-c_3) \gamma(r)
\end{split}
\end{equation}
where the last inequality holds if we can choose $\vep$ small enough so that the denominator in the third inequality is arbitrarily close to one. 
The constraint \eqref{e:betadef} on $\beta$ guarantees that we can do this.

Gathering together the bounds \eqref{p:split1}, \eqref{p:split2}, \eqref{p:bigedges}, \eqref{p:b1bd}, and \eqref{p:b2bd} gives the bound in the proposition's hypothesis for $\vep > r^{\xi-\rho}$.
\qed
\medskip

All that remains is to prove Claim \ref{claim}. We only give an outline of the proof of Claim \ref{claim} here, since, as mentioned before, the proof of the claim is very similar to the proof of \cite[Theorem~2]{KozNac11}. In the appendix we give all the necessary modifications to that proof so that it proves Claim \ref{claim}, but we omit many details.
\medskip

\appendix
\section{An outline of the proof of Claim \ref{claim}}\label{app}
This appendix contains an abridged version of the proof of \cite[Theorem 2]{KozNac11}. It is \emph{not} self-contained: it is intended to be read alongside \cite{KozNac11}. 
Indeed, the difference between the proof of \cite[Theorem 2]{KozNac11} and what we need is only superficial. The proof makes hardly any direct use of properties of percolation that do not also hold for LRP, and in the few places where it does, these properties are not crucial.

Let us start with some general remarks about the modifications.
The main idea of the proof of \cite[Theorem 2]{KozNac11} is the same as the proof of Claim \ref{claim:brw} above. Namely, we use the second moment method to bound the probability that there are many vertices on the boundary $\Xcal_{j}$ but that they have a combined cluster in $Q_{j+L} \setminus Q_j$ that is not too big. For BRW we could use the independence of the offspring of different particles to bound the first and second moment. For percolation this independence does not exist, because two vertices may for instance connect to the same cluster. To deal with these dependencies, Kozma and Nachmias perform a `regularity analysis' on the boundary: they show that the clusters emanating from most vertices on the boundary are close to independent (i.e., `regular') by showing that there are only a few vertices that are not regular. (Roughly speaking, a vertex in $\Ccal$ is regular if it has close to the expected number of other vertices in $\Ccal$ in it's vicinity, and if these vertices are `nicely spread out'.) For the proof of Claim \ref{claim} we need a similar analysis.

One of the first things one might notice when comparing this paper with \cite{KozNac11} is that \cite{KozNac11} frequently states events, summations, etc.\ in terms of $\partial Q_r$, the surface of the cube $Q_r$. In our setting $\partial Q_r$ does not mean very much, since LRP can produce very sparse clusters that may leave many surfaces $\partial Q_r$ empty while still reaching a great distance. For this reason we consider $Q_r^c$ instead when defining the one-arm event. In \cite{KozNac11} also uses the surface $\partial Q_j$ in regeneration arguments. Clearly that does not work for LRP either, since $X_j$ and $A_j$ are defined in terms of $\ann$. The problem is easily solved by replacing $\partial Q_j$ with $\ann$. 

Before we proceed with a detailed treatment, let us remark that these are the key modifications to the proof of \cite[Theorem 2]{KozNac11}:\begin{itemize}
	\item Replace $\partial Q_r$ with $Q_r^c$ in the definition of the one-arm event.
	\item Replace $\partial Q_j$ with $\ann$ when using $\partial Q_j$ as a regeneration surface.
	\item Replace the ersatz bound $\gamma(r) \le C /\sqrt{r}$ (which does not hold for LRP) with the bound from Lemma \ref{lem:subopt}, $\gamma(r) \le C/r^\xi$ with $\xi = \left( \frac{\alpha}{2\alpha + 2} \wedge \smallfrac{2}{5}\right)$ (which does hold). To accommodate this replacement we need to `stretch' some length scales by a factor~$1/\xi$.
	\item Replace the two-point function bound $\Ppc(x \conn y) \asymp \|x-y\|_2^{2-d}$ with the bound $\Ppc(x \conn y) \asymp \|x-y\|_2^{\twa-d}$ (i.e., \eqref{e:twopt}).
	\item Replace exponents that are integer multiples of $2$ with exponents that are multiples of $\twa$ if the objects in question are related to the volume of a cluster, and replace them with multiples of $\rho$ if the exponents are related to the one-arm event. (Powers of logarithms will remain unchanged.)
	\item Replace $x+Q_L$ by $H_{L,j}(x)$ (as defined in \eqref{e:shcube} below) where needed.
\end{itemize}

In what follows, we give the required modifications to the proof of \cite[Theorem 2]{KozNac11} per section (the proof spans \cite[Sections 3, 4, and 5]{KozNac11}). 
For ease of reference we will maintain the numbering scheme of \cite{KozNac11} in this appendix, but affix the letter `M' to the label to indicate that it is a modification.

\subsubsection*{Modifications for \cite[Section 3]{KozNac11}}
This section proves \cite[Lemma 1.1]{KozNac11}, which states a lower bound on a connection probability. This bound is not used directly. What is used is the lemma's corollary, \ck{Corollary 3.2}. This corollary is used only in the proof of \ck{Lemma 4.5}. 
It turns out that for LRP the bound that is needed in the proof of \ck{Lemma 4.5} is trivial, since connection probabilities in LRP can always be bounded from below by the probability of a connection by a single edge, which is large enough. Hence we may skip \ck{Section 3}.
	
\subsubsection*{Modifications for \cite[Section 4.1]{KozNac11}}
This section states the regularity result \ck{Theorem 4}, modified below. Roughly speaking, this theorem states that it is unlikely that the cluster has many `dense patches' on the boundary of a cube $Q_j$.
 
 We need to make several modifications in this section. We start by redefining `typical clusters':
\[
	\Tcal_s (x) :=\{|\Ccal(x) \cap (x + Q_s)| < s^{2 \twa} \log^7 s\}.
\]
Recall that $\Ccal(x; Q_{j})$ is defined as the connected cluster of the vertex $x$ in the modified configuration where all edges that do not have both vertices in $Q_j$ are closed.
With $\ann$ replacing $\partial Q_j$ and the new definition of $\Tcal_s(x)$, the definitions of $s$-bad and $K$-irregular are essentially the same:
\pagebreak
\begin{customdef}{M4.1} For $x \in  \ann$ and positive integers $s$ and $K$:
\begin{itemize}
	\item[(i)] We say that $x$ is $s$-bad if $\Ccal(x;Q_{j})$ satisfies
	\[
		\Ppc(\Tcal_s(x) \mid \Ccal(x; Q_{j})) \le 1- \e^{-\log^2 s}.
	\]
	\item[(ii)] We say that $x \in  \ann$ is $K$-irregular if there exists $s \ge K$ such that $x$ is $s$-bad.
\end{itemize}
\end{customdef}
We redefine the number of $K$-irregular vertices at level $j$,
\[
	X_j^{K\mathrm{-irr}} := \# \{x \in \ann : 0 \lstackrel{Q_{j}}{\conn} x \text{ and $x$ is $K$-irregular}\}.
\]
With these redefinitions we can leave \ck{Theorem 4} otherwise unchanged:
\begin{customthm}{M4}\label{thm:knreg} There exist constants $C, c > 0$ such that for any $K$ sufficiently large and any $j$ and $M$ we have
\emph{
\[
	\Ppc\left( X_j \ge M \text{ and }  X_{j}^{K\mathrm{-irr}} \ge \smallfrac12 X_j \right) \le C j^d \e^{-c \log^2 M}.
\]}
\end{customthm}

\subsubsection*{Modifications for \cite[Section 4.2]{KozNac11}}	
To verify whether a vertex $x$ is $K$-regular we need to inspect the entire percolation configuration. In this section, the notion of $K$-local-regularity is introduced, and it is shown that $K$-local-regularity implies $K$-regularity. The advantage of $K$-local-regularity is of course that we only need to inspect a finite region around $x$ to verify whether $x$ is $K$-locally-regular. We again need to make several modifications.

Recall from Lemma \ref{lem:subopt} that $\gamma(r) \le C/r^\xi$ with $\xi = \left( \frac{\alpha}{2\alpha + 2} \wedge \smallfrac{2}{5}\right)$.
We redefine what it means to be locally regular:
\begin{customdef}{M4.2} For $x \in \ann$ and a positive integer $s$ we say that the event $\Tcal_{s}^{\mathrm{loc}}(x)$ occurs if the following two events occur:
\begin{enumerate}
	\item For all $y \in x + Q_s$,
	\[
		|\Ccal(y ; x + Q_{s^{2d/\xi}}) \cap (x + Q_s)| < s^{2 \twa} \log^4 s.
	\]
	\item There exist at most $\log^3 s$ disjoint open paths starting in $x + Q_s$ and ending at $x + Q_{s^{2d/\xi}}^c$.
\end{enumerate}
\end{customdef}

\begin{customclaim}{M4.1}
For any $x \in \Zd$ and positive integer $s$,
\[
	\Tcal_{s}^{\mathrm{loc}}(x) \Rightarrow \Tcal_s (x).
\]
\end{customclaim}
The proof of \cite[Claim 4.1]{KozNac11} holds, mutatis mutandis.

We continue.
\begin{customdef}{M4.3} For $x \in  \ann$ and positive integers $s$ and $K$ we define:
\begin{itemize}
	\item[(i)] We say that a cluster $\Ccal$ in $B := (x + Q_{s^{4d^2/\xi}}) \cap Q_{j}$ is a `spanning cluster' if $x \in \Ccal$ or if $\Ccal \cap Q_{j}$ intersects both $x + Q_{s^{4d^2/\xi}}^c$ and $x+Q_{s^{2d/\xi}}$.
	\item[(ii)] We say that $x$ is $s$-locally-bad if there exist spanning clusters $\Ccal_1, \dots, \Ccal_m$ in $B$ such that
	\[
		\Ppc(\Tcal_s^{\mathrm{loc}}(x) \mid \Ccal_1, \dots, \Ccal_m) \le 1 - \e^{-\log^2 s}.
	\]
	\item[(iii)] We say that $x \in \ann$ is $K$-locally-irregular if there exists $s \ge K$ such that $x$ is $s$-locally-bad.
\end{itemize}
\end{customdef}
We modify \ck{Claim 4.2}:
\begin{customclaim}{M4.2} For any $x \in \ann$ and a positive $s$ we have that if $x$ is $s$-locally-good, then $x$ is $s$-good.
\end{customclaim}
Claim M4.2 has the same proof as \ck{Claim 4.2}, mutatis mutandis.
\medskip

We redefine
\[
	X_j^{K-\mathrm{loc-irr}} := \# \{x \in \ann : 0 \lstackrel{Q_{j}}{\conn} x \text{ and $x$ is $K$-locally-irregular}\}.
\]
With these new definitions, we can leave \ck{Theorem 5} otherwise unchanged:
\begin{customthm}{M5}
There exist constants $C > c > 0$ such that for any $K$ sufficiently large and any $j$ and $M$ we have
\[
	\Ppc \left(X_j \ge M \text{ and } X_j^{K-\mathrm{loc-irr}} \ge \smallfrac12 X_j \right) \le C j^d \e^{-c \log^2 M}.
\]
\end{customthm}
Theorem M4 follows directly from Theorem M5 and Claim M4.2.

\subsubsection*{Modifications for \cite[Section 4.3]{KozNac11}} This section states two lemmas that give a large deviations type bound on the probability that a vertex is $s$-locally bad.
We state two modified lemmas:
\begin{customlem}{M4.3} For $x \in \ann$ and positive integer $s$ we have
\emph{
\[
	\Ppc(x \text{ is $s$-locally-bad}) \le C \e^{-c \log^4 s}.
\]}
\end{customlem}
\begin{customlem}{M4.4} There exists some constant $c>0$ such that for all $s>0$ and $\lambda>0$ we have
\[
	\Ppc(\max_{y \in Q_s} | \Ccal(y) \cap Q_s| > \lambda s^{2 \twa}) \le s^{d-3 \twa} \e^{-c \lambda}.
\]
\end{customlem}
Lemma M4.4 is used to prove Lemma M4.3.
Note that the exponents $4$ and $6$ in \cite[Lemma 4.4]{KozNac11} have been replaced here with $2\twa$ and $3\twa$, corresponding to the fact that this is a volume estimate.
Write $\Ccal_{\mathrm{max}}$ for the maximum in Lemma M4.4. The proof of Lemma M4.4 is the same as \ck{Lemma 4.4}, except that we now use the bound
\[
	\Epc[\Ccal_{\mathrm{max}}^k] \le k! C_1^k s^{d - 3\twa + 2\twa k},
\]
which follows from \cite[\S 4.3, Lemma 2]{Aize97} (Aizenman's $\eta$ is $\twa - 2$ in our case, corresponding to the assumption \eqref{e:twopt}).
\medskip

The proof of \ck{Lemma 4.3} does not work for Lemma M4.3 because it uses the bound $\gamma(r) \le \Ppc(|\Ccal| \ge r/\Lambda) \le C/\sqrt{r}$, which does not hold for LRP. Therefore we need to redo two of the bounds using Lemma \ref{lem:subopt}, but the rest of the proof holds, mutatis mutandis.

The first of these bounds is used to bound the probability of the complement of the event in Definition M4.2(b): We use \eqref{e:subopt}, translation invariance, monotonicity of the one-arm event, and the union bound to bound
\[
	\Ppc(x + Q_s \conn x+Q_{s^{2d/\xi}}^c) \le Cs^d \Pp(0 \conn Q_{s^{2d/\xi}-s}^c) \le \frac{Cs^d}{(s^{2d/\xi}-s)^\xi} \le \frac{C}{s}.
\]
The second bound goes along exactly the same lines: the probability that there exist at least $\log^3 s$ spanning clusters is at most
\[
	\left(Cs^{2d/\xi}\right)^{\log^3 s} \left(C (s^{4d^2/\xi} - s^{2d/\xi})^{-\xi} \right)^{\log^3 s} \le C \e^{-c \log^4 s}.
\]
All other parts of the proof are the same.

\subsubsection*{Modifications for \cite[Section 4.4]{KozNac11}} This section presents an exploration process that is used to define two martingales. These martingales are used to control the number of vertices and $s$-locally-bad vertices on the boundary of $Q_j$. To define the exploration for LRP we need a different notion of a connection event `through' a set $B$:

For a vertex $x$ and vertex sets $A$ and $B$, we write $\{x \lstackrel{B}{\lrsq} A\}$ for the event that there exists a path of open edges from $x$ to $A$ such that all vertices along that path, except possibly the terminal vertex of the path in $A$, belong to $B$. 
Another way of looking at this event for the case $A \cap B = \emptyset$, is that there exists a path from $x$ to some $y \in B$ in the modified percolation configuration where all edges that do not have both end-points in $B$ are are made closed, and there exists an open edge from $y$ to a vertex $z \in A$ in the unmodified configuration.

We redefine the set of boxes
\[
	G(w) := \{(Q_{2s^{4d^2/\xi}} + v) \cap Q_{j} \, : \, v \in (4 s^{4 d^2/\xi}+1) \Zd + w\} \setminus \emptyset,
\]
and the sets of explored and active boxes:
at time $1$,
\[
\begin{aligned}
	E_1 &:= \{q \in G(w) : q \cap Q_{j-\delta r}^c = \emptyset\},\\
	A_1 &:= \{q \in G(w) : \exists x \in q \text{ s.t. } 0 \lstackrel{\cup E_1}{\lrsq} x\} \setminus E_1,
\end{aligned}
\]
and at time $i$:
\[
\begin{aligned}
	E_i &:= E_{i-1} \cup \{q_i\},\\
	A_i &:= \left( A_{i-1} \cup \{q \in G(w) : \exists x \in q \text{ s.t. } 0 \lstackrel{\cup E_i}{\lrsq} x\}\right) \setminus E_i.
\end{aligned}
\]
We say that $q \in G(w)$ is $s$-bad if there exists an $x \in \ann$ that is $s$-locally bad and such that $(x + Q_{s^{4d^2/\xi}}) \cap Q_{j} \subset q$. The definition of the filtration $\{\mathcal{F}_i\}$ and the martingale $(\beta_i)$ remain unchanged: The filtration $\{\Fcal_i\}$ is the configuration restricted to $\bigcup_{j \le i}E_j$, and
\[
	\beta_i = \beta_{i-1} + \indi_{\{q_i \text{ is $s$-bad}\}} - \Ppc(q_i \text{ is $s$-bad} \mid \Fcal_{i-1}),
\]
with $\beta_1 =0$.
We redefine the martingale $(\gamma_i)$: start with $ \gamma_1 =0$ and let
\[
	\gamma_i = \gamma_{i-1} + \indi_{\left\{\exists x \in q_i \cap \ann \,:\, 0 \lstackrel{\cup E_i}{\lrsq} x\right\}} - \Ppc(\exists x \in q_i \cap \ann \,:\, 0 \lstackrel{\cup E_i}{\lrsq} x \mid \Fcal_{i-1}).
\]
We restate \ck{Lemma 4.5} and its proof here with the needed modifications:
\begin{customlem}{M4.5} There exist constants $C_1 > 0$ and $c_1 > 0$ such that for any $j$, $s$, and $M$ we have
\[
	\Ppc(c_1 \e^{-C_1 \log^2 s} \tau \ge X_j \ge M) \le C \e^{-c M + C \log^2 s},
\]
where $\tau$ is the stopping time of the exploration defined above.
\end{customlem}
\proof[Proof of Lemma M4.5] For every $i$,
\[
	X_j \ge \gamma_i + \sum_{k=1}^i \Ppc(\exists x \in q_k \cap \ann  \,:\, 0 \lstackrel{\cup E_k}{\lrsq} x \mid \Fcal_{k-1}).
\]
Since we explored $q_k$, there exists some $z \in q_k^c$ such that $\{0 \lstackrel{\cup E_{k-1}}{\conn} z\}$ occurs and $\{z \lstackrel{q_k}{\lrsq}x\}$ occurs, so given $\Fcal_{k-1}$, the probability that there exists $x \in q_k \cap \ann$ such that $\{0 \lstackrel{\cup E_k}{\lrsq} x\}$ occurs is at least the probability of $\{z \lstackrel{q_k}{\lrsq} x\}$. The probability of this event is bounded from below by the probability that there exists an open edge between $z$ and $x$, which, by Definition \ref{def:D}(c), is at most $c (s^{4d^2/\xi})^{-d-\alpha} \ge c_2 \e^{-C_2 \log^2 s}$. From here the proof continues as in \cite{KozNac11}. Observe that this short argument completely circumvents the calculations in \cite[Section 3]{KozNac11} (but it only works for LRP).\qed
\medskip

Redefining
\[
	X_j^{s\textrm{-loc-bad}} := \#\{x \in \ann \, : \, x \text{ is $s$-locally-bad}\}
\]
and
\[
	X_j^{s\textrm{-loc-bad}}(w) := \# \{x \in \ann \, : \, x \text{ is $s$-locally-bad, and }\exists q \in G(w) \text{ s.t. } x \in q\},
\]
the statement and proof of \cite[Lemma 4.6]{KozNac11} becomes
\begin{customlem}{M4.6} There exist constants $C_3 >0$ and $c_3>0$ such that for any $j,$ $s$, $w$, and $M$, and for any real number $\mu \ge C_3 \e^{-c_3 \log^4 s}$ we have
\emph{
\[
	\Ppc( \mu^{-1} X_{j}^{s\textrm{-loc-bad}}(w) \ge \tau \ge M) \le C s^{8d^2/\xi}\mu^{-2} \exp\left(-c s^{-8 d^2 \xi} \mu^2 M \right).
\]}
\end{customlem} 
The proof is basically the same as given in \cite{KozNac11}, but the proof given there contains a small mistake, so we give a short version of the proof to indicate how the different exponents come about.

\proof By Lemma M4.3 it follows that
\[
	\Ppc(q_k \text{ is bad} \mid \Fcal_{k-1}) \le \sum_{x \,:\, x+Q_{s^{4d^2 /\xi}} \subset q_k} C \e^{-c \log^4 s} \le \smallfrac12 C_3 \e^{-c_3 \log^4 s}.
\]
It thus holds deterministically that
\[
	X_j^{s\textrm{-loc-bad}}(w) \le (2s^{4d^2/\xi}+1)^d |\{k \le \tau \, : \, q_k \text{ is bad}\}| \le (2s^{4d^2/\xi}+1)^d \left(\beta_\tau + \smallfrac12 C_3 \e^{-c_3 \log^4 s} \tau \right)
\]
(In \cite{KozNac11} this bound is stated with a factor $s^{d-1}$ multiplying only the first term, but this cannot hold deterministically even in the nearest neighbor case. It does not matter: the polynomial factor is negligible compared to the exponential factor.) Now $X_j^{s\text{-loc-bad}}(w) \ge \mu \tau$ implies
\[
	\beta_\tau \ge c s^{-4d^3 /\xi} \tau \left(\mu - \smallfrac12 C_3 \e^{-c_3 \log^4 s} \right) \ge c_4 s^{-4d^3/\xi} \mu \tau.
\]
Continuing the proof with this bound we find
\[\begin{split}
	\Ppc(\mu^{-1} X_j^{s\text{-loc-bad}}(w) \ge \tau \ge M) &\le \Ppc\left(c_4^{-1} s^{4d^3/\xi} \mu^{-1} \beta_\tau \ge \tau \ge M \right)\\
	& \le \sum_{i=M}^\infty \Ppc(\beta_i \ge c_4 s^{-4d^3/\xi} \mu i)\\
	& \le \sum_{i=M}^\infty \exp\left(-c s^{-8d^3 /\xi} \mu^2 i \right)\\
	& \le C s^{8d^2/\xi}\mu^{-2} \exp\left(-c s^{-8 d^2 \xi} \mu^2 M \right). \qed
\end{split}\] 

Theorem M5 is proved using Lemmas M4.5 and M4.6, and goes as in \cite{KozNac11}, mutatis mutandis.

\subsubsection*{Modifications for \cite[Section 5]{KozNac11}}
This section uses a second moment method to bound the probability that the percolation cluster contains many vertices on the boundary of $Q_j$, given that there are precisely $M$ $K$-regular vertices. This is combined with Theorem~M4 (which bounds the probability that there are too many $K$-irregular vertices) to give the desired bound.
The only changes needed in this section pertain to the volume and one-arm exponents (i.e., replacing $2$'s with $\twa$'s or $\rho$'s).

We start by defining the \emph{shifted outward-facing half-cube}
\begin{multline}\label{e:shcube}
	H_{L,j}(x) : = \{y \in \Zd : \|x + \sign(x_{\iota_x})(j-\|x\|_\infty) e_{\iota_x} -y\|_\infty \le L^{\rho/\twa}\\
	 \text{ and } |y_{\iota_x}| > j \text{ and } \sign(y_{\iota_x}) = \sign({x_{\iota_x}})\},
\end{multline}
where $e_{\iota_x}$ is the unit vector in direction $\iota_x$  (similar to \eqref{e:hcube}). Observe that $H_{L,j}(x) \subset Q_{j+L} \setminus Q_j$ for all $x \in \ann$ (this is the result of the extra outward shift with respect to \eqref{e:hcube}).

We redefine the notion of $(j,L,K)$-admissible to be a pair of vertices $(x,y)$ that satisfies:
\begin{itemize}
	\item $x \in \ann$ and $y \in H_{L,j}(x)$.
	\item $0 \lstackrel{Q_j}{\conn} x$ and $x \conn y$.
	\item $x$ is $K$-regular.
	\item There exists a vertex $\tilde x \in Q_j^c$ such that $\{x,\tilde x\}$ is an open and pivotal edge for the event $0 \conn y$. (If there is more than one such edge, choose the lexicographically first one).
\end{itemize}

With $Y(j,L,K)$ denoting the number of pairs $(x,y)$ that are $(j,L,K)$-admissible, \cite[Lemmas 5.1 and 5.2]{KozNac11} can now be restated as follows:
\begin{customlem}{M5.1} Let $K$ be sufficiently large, and let $j$, $M$ and $L$ be integers such that $M \ge L^\rho/2$. Then there exists a constant $c = c(K) >0$ such that
\emph{
\[
	\Epc [Y(j,L,K) \indi_{\{X_j^{K\text{-reg}} = M\}}] \ge c M L^\rho \Ppc(X_j^{K\text{-reg}} = M).
\]}
\end{customlem}
\begin{customlem}{M5.2} Let $j$, $K$, $M$ and $L$ be integers. Then
\emph{
\[
	\Epc [Y(j,L,K)^2 \indi_{\{X_j^{K\text{-reg}} = M\}}] \le C M^2 L^{2\rho} \Ppc(X_j^{K\text{-reg}} = M).
\]}
\end{customlem}

From here on the proof of Claim \ref{claim} easily follows:
Using Theorem M4 we bound
\begin{equation}\label{e:thm2pf}
\begin{split}
	\Ppc(X_j \ge L^\rho \text{ and } A_j \le aL^{2\rho}) &\le \Ppc(X_j \ge L^\rho \text{ and } X_j^{K\text{-irr}} \ge \smallfrac12 L^{\rho})\\
	& \le C j^d \e^{-c \log^2 j} + \sum_{M \ge L^\rho/2} \Ppc(X_j^{K\text{-reg}} = M \text{ and } A_j \le c L^{2 \rho})
\end{split}\end{equation}
(here we also used that $L \ge j^\theta$ for some $\theta >0$.)

Now we use that $A_j \ge Y(j,L,K)$ by definition,\footnote{
Note that there is a small mistake in \cite{KozNac11}: in \cite[(5.2)]{KozNac11} it is claimed that $A \ge Y$. Using their definitions this inequality does not hold since the random variable $Y(j,L,K)$ may count vertices outside of $Q_{j+L} \setminus Q_j$, while $A_j$ does not. Our use of $H_{L,j}(x)$ rather than $x + Q_L$ in the definition of $(j,L,K)$-admissible takes care of this.} and we use the Paley-Zygmund inequality and Lemmas M5.1 and M5.2 to bound
\[
		\Ppc(X_j^{K\text{-reg}} = M \text{ and } A_j \le c L^{2 \rho}) \le (1-c')\Ppc(X_j^{K\text{-reg}} = M).
\]
Substituting this bound in \eqref{e:thm2pf}, using that the exponential term is negligible and using that $\{X_j^{K\text{-reg}} \ge \smallfrac12 L^{\rho}\}$ implies $\{0 \conn Q_r^c\}$, we complete the proof of Claim \ref{claim}.
\medskip

It remains to prove Lemmas M5.1 and M5.2. These proofs can be given exactly as in \cite{KozNac11} if we consistently adhere to the following substitutions:
\begin{itemize}
	\item As above, we should always use $x \in \ann$ instead of $x \in \partial Q_j$ and $y \in H_{L,j}(x)$ instead of $y \in x + Q_L$.
	\item Two-point function estimates of the form $\|x-y\|_2^{2-d}$ should be replaced with $\|x~-~y\|_2^{\twa-d}$.
	\item Corresponding to the above two substitutions, all exponents pertaining to volume estimates should be replaced with their corresponding multiples of $\rho$ if they arise due to summing over  sets of the form $H_{L,j}(x)$, and with $\twa$ if they are due to bounds that use $\Tcal_s$.
\end{itemize}
For instance, the estimate \cite[(5.5)]{KozNac11} becomes
\[
	\sum_{y \in H_{L,j}(x)} \Ppc(x' \conn y \text{ only on }A) \le C L^\rho \sum_{z \in A} \|z - x'\|_2^{\twa -d},
\]
and the claim \cite[(5.7)]{KozNac11} becomes
\[
	|A_t| < 2^{2\twa (t+1)} (t+1)^7 \qquad \forall t \text{ such that } 2^t \ge K/2.
\]

Indeed, in these final steps we estimate the first and second moment of the number of $(j,L,K)$-admissible pairs by using the properties of $K$-regularity (i.e., that it is unlikely that there are many vertices on the boundary near a $K$-regular vertex).

As a last remark: at the end of the proof of Lemma M5.1 we should note that the local modification involves not $(2K)^d$ as in the proof of \ck{Lemma 5.1}, but $(2K)^{2d}$ edges. This is of course of no importance.

This concludes the description of all the ways that the proof of \cite[Theorem 2]{KozNac11} should be modified to prove Claim \ref{claim} above.

\subsection*{Acknowledgments:} The author wishes to thank Remco van der Hofstad and J\'{u}lia Komj\'{a}thy for stimulating discussions and many useful comments, and Asaf Nachmias and Akira Sakai for their useful remarks.

\begin{small}
\bibliographystyle{abbrv}
\bibliography{../bib/TimsBib}
\end{small}
\end{document}